\documentclass[preprint,11pt]{elsarticle}

 \usepackage{amsmath,amsthm,amssymb,enumerate}
\usepackage{gensymb}

\usepackage[sc]{mathpazo}
\usepackage{multirow} 

\usepackage{multicol}

\usepackage{newtxtext,newtxmath}
\usepackage{subfig}
\usepackage{graphicx}
\usepackage{epstopdf}
\usepackage{hyperref}
\usepackage{cancel}
\usepackage[margin=3cm]{geometry}
\usepackage{tikz}
\usepackage{caption}
\usetikzlibrary{shapes,calc}
\usepackage{verbatim}
\usepackage{mathrsfs}
\usepackage{algorithm}
\usepackage{accents}

\newtheorem{theorem}{Theorem}[section]
\newtheorem{lemma}[theorem]{Lemma}
\newtheorem{proof of lemma}[theorem]{Proof of Lemma}

\theoremstyle{definition}

\newtheorem{example}[theorem]{Example}

\newtheorem{remark}[theorem]{Remark}
\numberwithin{equation}{section}

\def \R{{{\Bbb R}}}

\allowdisplaybreaks

\def\R{\mathbb{R}}

\def\dx{{\rm~dx}}

\def\u{\boldsymbol{u}}
\def\V{\mathbf{V}}
\def\Ul{\underline{\boldsymbol{U}}}
\def\bvl{\underline{\boldsymbol{v}}}
\def\ul{\underline{\boldsymbol{u}}}
\def\zl{\underline{\boldsymbol{z}}}
\def\Il{\underline{\boldsymbol{I}}}
\def\wl{\underline{\boldsymbol{w}}}
\def\w{\boldsymbol{w}}
\def\x{\boldsymbol{x}}
\def\G{\boldsymbol{G}}
\def\r{\boldsymbol{r}}

\def\f{\boldsymbol{f}}

\def\div{{\rm div~}}
\def\bv{\boldsymbol{v}}
\def\z{\boldsymbol{z}}
\def\b{\boldsymbol{b}}

\def\bn{\boldsymbol{n}}

\def\V{\boldsymbol{V}}

\def\cM{\mathcal{M}}
\def\cT{\mathcal{T}}
\def\cF{\mathcal{F}}

\def\jump#1{\left[\hskip -3.5pt\left[#1\right]\hskip -3.5pt\right]}

\newcommand{\norm}[1]{\left\Vert #1 \right\Vert}

\newcommand{\tnorm}[1]{\ensuremath{\left| \! \left| \! \left|} #1\ensuremath{\right| \! \right| \! \right|}}
\newcommand{\tnormn}[1]{\ensuremath{|\!|\!|} #1\ensuremath{|\!|\!|}}

\title{A Local Projection Stabilised HHO Method for the Oseen Problem}

\begin{document}
\begin{frontmatter}
\author{Gouranga Mallik\footnote{Department of Mathematics,  School of Advanced Sciences,\\
Vellore Institute of Technology, Vellore - 632014; gouranga.mallik@vit.ac.in}, 
Rahul Biswas\footnote{Department of Mathematics, Indian Institute of Science, Bangalore - 560012; rahulbiswas@iisc.ac.in},
and Thirupathi Gudi\footnote{Department of Mathematics, Indian Institute of Science, Bangalore - 560012; gudi@iisc.ac.in}}
	
\begin{abstract}
In this article, we consider a local projection stabilisation for a Hybrid High-Order (HHO) approximation of the Oseen problem. We prove an existence-uniqueness result under a stronger SUPG-like norm.
We improve the stability and provide error estimation in stronger norm for convection dominated Oseen problem. We also derive an optimal order error estimate under the SUPG-like norm for equal-order polynomial discretisation of velocity and pressure spaces. Numerical experiments are performed to validate the theoretical results.
\end{abstract}

\end{frontmatter}


\section*{Introduction}\label{sec:Intro}
The Navier-Stokes equation models the flow of fluid in a domain. A solution to these equations is important in many engineering problems. Linearizing and time-discretizing the Navier-Stokes equation, we obtain the Oseen problem:
\begin{equation}\label{Oseen_eqn_intro}
\begin{aligned}
-\epsilon \Delta \u + (\b \cdot \nabla) \u + \sigma \u  + \nabla p =&\f \text{ in } \Omega, \\
\div \u =&0 \text{ in } \Omega\\
\u =& 0 \text{ on } \partial \Omega,
\end{aligned}
\end{equation}
where $\u$ denotes the velocity of the fluid and $p$ denotes the pressure. Here $\Omega \subset \R^d$ with $d\in\{2,3\}$ is a bounded polytopal domain with Lipschitz boundary $\partial \Omega$. The force function $\f$ is in $[L^2(\Omega)]^d$. The viscosity coefficient is denoted by $\epsilon$, where $0 < \epsilon \ll 1$. The convection coefficient $\b$ is a $[W^1_\infty(\Omega)]^d$ function such that $\div \b=0$. The reaction coefficient is a positive constant denoted by $\sigma$.

Fluid flow problems with dominant convection produce boundary and interior layers. It is well-known that the numerical solution to these problems using the usual Galerkin method cannot capture these small layers. Instead, they produce nonphysical solutions that contain spurious oscillations. To eliminate the effect of convection, one can add stabilisations. The Hybrid High-Order method with upwind stabilisation of \cite{Oseen_HHO} eradicates the oscillation. In this article, we focus on local projection stabilisation using a Hybrid High-Order approximation on general polygonal meshes and establish the stability and error estimation in a stronger norm.

In the last few years, there has been a growing interest in high-order polynomial approximations of solutions to PDEs on general polytopal meshes. Due to the vast literature in this area, we cite a few well-known works: the Hybridizable Discontinuous Galerkin (HDG) method in \cite{Cock_Gopal_Laz_09_Unif,Piet_Ern_12_DG_book}. The HDG has been extended further to the convection diffusion problem in \cite{Cock_Dong_09_HybridDG} and the Oseen problem in \cite{Oseen_HDG_cock}. The Virtual Element Method (VEM) has been studied in \cite{Beira_Brezz_Cang_Man_Mar_13_VEM,Beira_Brezz_Mari_13_VEM_Elast,Brezz_Falk_Mari_14_VEM_Mixed}. The VEM has also been applied to the convection diffusion problem with SUPG stabilisation \cite{VEM_con_diff_supg} and to the Oseen problem with LPS stabilisation \cite{Oseen_VEM_LPS}. The Weak Galerkin method is introduced in \cite{Mu_Wang_Ye_15_WeakGaler_Nonlin,Wang_Ye_13_Weak_Galerkin,Wang_Ye_14_Mixed_Weak_Galerkin} and the Gradient Discretisation method is introduced in \cite{Dron_Eym_Herb_GDM_16,Dron_Eym_GDM_18_book,Piet_Dron_Man_18_GradDis_Polytope}. The Multiscale Hybrid-Mixed method has been studied in \cite{Arya_Hard_Pare_Vale_13_Multiscale}. The focus of our article is on the Hybrid High-Order (HHO) method, originally introduced in \cite{Piet_Ern_Lem_14_arb_local,ern_pietro_elasticity}. For an overview of the HHO method, we refer to \cite{Piet_Jero_HHO_Book_20}. HHO is a robust method based on local polynomial reconstruction. It is independent of the dimension of the problem and suitable for local static condensation, which drastically reduces the computational cost of the matrix solver.

The HHO method is closely related to the HDG method, but differs in the choice of stabilisation; see \cite{Cock_Piet_Ern_16_HHO_HDG} for details. In the nonconforming Virtual Element methods (ncVEM) one takes the projection of the virtual function in the stabilisation, whereas the HHO method takes a reconstruction of the function in the stabilisation. In \cite{Lem_21_Bridging_HHO_VEM}, the connection of the HHO method with the virtual element method is discussed. See \cite{Bone_Ern_14_poly, Dron_Eym_Gall_10_unified,Brezzi_Lip_Shas_05_Conv_Mimetic, Brezzi_Lip_Shas_07_polyhedral, Kuzn_Lip_Shas_04_Mimetic_poly} for related works. In the lowest-order case ($k=0$), the HHO method resembles the Hybrid Mixed Mimetic family, hence the mixed-hybrid Mimetic Finite Differences, the Hybrid Finite Volume and the Mixed Finite Volume methods too, see \cite{Dron_Eym_Gall_10_unified, Brez_Lipn_Simo_05_Mim_Dif,Eym_Gall_Herb_10_NonconfMesh,Dron_14_FV_Diff, Dron_Eym_06_MixedFV}. We state the following predominant works on the HHO methods: problem with pure diffusion \cite{Piet_Ern_Lem_14_arb_local}, advection-diffusion problem \cite{Piet_Dron_Ern_15_HHO_Adv_Diff_Rea},  interface problems \cite{Burman_Ern_18_HHO_Interface}, for linear PDEs, elliptic obstacle problem \cite{Cicu_Ern_Gudi_20_HHO_Obstacle}, Stokes problem \cite{Piet_Ern_Link_16_HHO_Stokes}, the Oseen problem \cite{Oseen_HHO}, and the steady incompressible Navier-Stokes equations~\cite{Piet_Krell_18_NSE}.

In the study of stabilisations for fluid flow problems with a high Reynolds number, the SUPG method by Hughes and Brooks \cite{MR679322} is the most well known. SUPG is studied for the incompressible Navier-Stokes equations in \cite{MR1377246}. There is a wide range of stabilisation techniques in the literature, some of them are: the least squares method \cite{MR2454024}, residual free bubbles technique \cite{MR2454024}, continuous interior penalty method \cite{MR2454024}, and the discontinuous Galerkin method \cite{MR2454024}. In this article, we are interested in the local projection stabilisation scheme, originally introduced for the Stokes problem by Becker and Braack \cite{MR1890352} for the Stokes problem. It has also been studied for the transport problem, the scalar convection diffusion problem, the Oseen problem, and the Navier-Stokes equations \cite{MR2121360,MR2454024,BRAACK2006372,MR2206447,MR2242376}. A nonconforming patchwise LPS method using Crouzeix-Raviart elements for the convection diffusion problem has been studied in \cite{MR3915226} and for the Oseen problem in \cite{biswas2019edge}. Knobloch has studied a generalised version of LPS for the convection diffusion problem in \cite{MR2670000} and also for the Oseen problem \cite{MR3003108}. 

The SUPG method naturally gives an additional control over the advective derivative of the velocity; however, the usual LPS methods in \cite{MR2121360,MR2206447} do not provide this. Moreover, the LPS methods in \cite{MR2121360,MR2206447} work through a two-level mesh approach or through enrichment. The articles \cite{MR2121360,MR2206447,MR2670000,MR3003108} need to satisfy a local inf-sup condition necessary for error analysis and stability. In this article, we employ a generalised LPS technique  
to design an HHO method for the Oseen problem motivated by the works in \cite{Oseen_HHO} on the HHO approximation and in \cite{MR3003108} on LPS stabilization for the Oseen problem. The additional LPS term provides control on the advective derivative. Moreover, we employ a one-level approach that does not require any enrichment of discrete spaces and the need for a local inf-sup condition as seen in \cite{MR3003108}.

In this article, along with the LPS method, we have added another velocity stabilisation to control the normal jump of the solution. This helps to further stabilise the solution. Moreover, we also need pressure stabilisation to stabilise the pressure gradient. In a nutshell, the LPS-HHO method is a combination of a usual HHO method for the Oseen problem combined with the above stabilisations along with an upwind term. Compared to \cite{Oseen_HHO}, we have proven that the LPS stabilisation term in the formulation helps to prove a stability result under a stronger SUPG-like norm. Moreover, in this article, the presence of normal jump stabilisation in the discrete scheme gives epsilon robust error bounds. This can be seen in inequality \eqref{apriori_estimate_eqn_15}.


The rest of the article is organised as follows: Section~\ref{preliminaries} defines the Oseen problem along with some notation and preliminaries. Section~\ref{Discrete_Oseen_problem} deals with the discrete HHO formulation of the Oseen equations. Section~\ref{Wellposedness} provides the proof for the discrete well-posedness of the system in Section~\ref{Discrete_Oseen_problem}. Section 5 provides a priori error estimates. The numerical results are provided in Section 6. The notation $a\lesssim b$ means that there exists a generic constant $C$ independent of the meshsize such that $a\leq C b$. We abbreviate $a\lesssim b\lesssim a$ by $a\sim b$. For $M\subset \Omega$, the $L^2$-inner product on $L^2(M)$ is denoted by $(\cdot,\cdot)_{M}$ and $L^2$-norm by $\|\cdot\|_{M}$. We omit the subscript for the domain specification when $M=\Omega$. We extend the definition $(\cdot,\cdot)_M$ naturally to vector-valued functions as the sum of component-wise inner products. The analysis is done on standard $k$th order Sobolev spaces $H^k(\Omega)$ with the standard norm $\norm{\cdot}_k=\sum_{|\alpha|\leq k}\norm{D^{\alpha}\cdot}_{L^2(\Omega)}$. The Sobolev space $H^1(\Omega)$ with zero trace is denoted by $H^1_0(\Omega)$. The space of square-integrable functions with zero mean is denoted by $L^2_0(\Omega)$. For $M \subset \Omega, k\geq 0$ and $1\leq p\leq\infty$, let $\|.\|_{k,p,M}$ be the norm on the $k$th order Sobolev space $W^{k,p}(M)$. For $p=2$ we denote the norm by $\|.\|_{k,M}$.

\section{Continuous Problem, Notations and Preliminaries}\label{preliminaries}
In this section, we introduce the weak formulation of the Oseen problem~\eqref{Oseen_eqn_intro} and some preliminaries. Let $\V=[H^1_0(\Omega)]^d$ be the velocity space and $Q=L^2_0(\Omega)$ be the pressure space. The weak formulation for the Oseen problem~\eqref{Oseen_eqn_intro} is given by: Find $\u \in \V $ and $p \in Q$ such that 
\begin{equation}\label{Oseen_eqn_sep}
\begin{aligned}
a(\u,\bv)-b(p,\bv)&=(\f,\bv)\quad {\rm for ~all}~\bv\in \V,\\
b(q,\u)&=0 \quad {\rm for ~all}~q\in Q,
\end{aligned}
\end{equation}
where, the bilinear forms $a(\cdot,\cdot)$ and $b(\cdot,\cdot)$ are defined as
\begin{align*}
a(\u,\bv)&:= \epsilon(\nabla\u,\nabla\bv)+((\b\cdot \nabla) \u,\bv)+ (\sigma \u,\bv ),\\
b(q,\u)&:= (q,\div \u).
\end{align*}
Using the fact that $\div \b=0$ and $\sigma > 0$, one can show that the bilinear form $a(\u,\bv)$ is coercive. It is well known that the bilinear form $b(q,\u)$ is inf-sup stable for $\u \in [H^1_0(\Omega)]^d$ and $q \in L^2_0(\Omega)$. Therefore, the existence and uniqueness of the problem \eqref{Oseen_eqn_sep} can be shown using the Babu\v{s}ka-Brezzi condition, see \cite[Chapter IV]{MR851383}. An equivalent formulation for \eqref{Oseen_eqn_sep} seeks $(\u,p) \in \V \times Q$ such that 
\begin{align}\label{Oseen_eqn_weak}
A((\u,p),(\bv,q))=(\f,\bv), \quad \text{for all} \quad (\bv,q) \in \V \times Q,
\end{align}
where, the combined bilinear form is defined by
\begin{align*}
A((\u,p),(\bv,q)):= a(\u,\bv)-b(p,\bv)+ b(q,\u).
\end{align*}
The existence and uniqueness of the problem \eqref{Oseen_eqn_weak} can be proved in a similar manner. Henceforth, we will use this combined mixed formulation in our analysis.

\medskip

Consider a decomposition $\cT_h$ of the domain $\Omega$ that consists of a finite collection of nonempty disjoint open polyhedral cells $T$ such that $\overline{\Omega}=\cup_{T\in\cT_h} \overline{T}$ and $h=\max_{T\in\cT_h} h_T$, where $h_T$ is the diameter of $T$. We follow \cite{Piet_Jero_HHO_Book_20} for various notation and definitions related to the mesh sequence $\cT_h$. A closed subset $F$ of $\Omega$ is defined to be a mesh face if it is a subset of an affine hyperplane $H_F$ with positive $(d-1)$-dimensional Hausdorff measure and if either of the following two statements holds true: (i) There exist $T_1(F)$ and $T_2(F)$ in $\cT_h$ such that $F\subset\partial T_1(F)\cap \partial T_2(F)\cap H_F$; in this case, the face $F$ is called an internal face; (ii) There exists $T(F)\in\cT_h$ such that $F \subset \partial T(F) \cap \partial\Omega \cap H_F$; in this case, the face $F$ is called a boundary face. The set of mesh faces is a partition of the mesh skeleton, that is, $\cup_{T\in\cT_h}\partial T=\cup_{F\in\cF_h}\bar{F}$, where $\cF_h:=\cF_h^i\cup \cF_h^b$ is the collection of all faces, which is the union of the set of all internal faces $\cF_h^i$ and the set of all boundary faces $\cF_h^b$. Let $h_F$ denote the diameter of $F\in \cF_h$. For each $T\in \cT_h$, the set $\cF_T:=\{F\in \cF_h\, |\, F\subset\partial T\}$ denotes the collection of all faces contained in $\partial T$. Following \cite[Definition 1]{Piet_Ern_Lem_14_arb_local}, we assume that the mesh sequence $(\cT_h)_{h>0}$ is admissible in the sense that, for all $h>0$, $\cT_h$ admits a matching simplicial submesh $\mathfrak{T}_h$ (i.e., every cell and face of $\mathfrak{T}_h$ is a subset of a cell and a face of $\cT_h$, respectively) so that the mesh sequence $(\cT_h)_{h>0}$ is shape-regular in the usual sense and all the cells and faces of $\cT_h$ have a uniformly comparable diameter to the cell and face of $\cT_h$ to which they belong. Owing to \cite[Lemma~1.42]{Piet_Ern_12_DG_book}, for $T\in\mathcal{T}_h$ and $F\in\mathcal{F}_T$, $h_F$ is comparable to $h_T$ in the sense that
\begin{align*}
	\varrho^{2}h_T\leq h_F\leq h_T,
\end{align*}
where $\varrho$ is the mesh regularity parameter. Moreover, there exists an integer $N_{\partial}$ depending on $\varrho$ and $d$ such that (see \cite[Lemma~1.41]{Piet_Ern_12_DG_book})
\begin{align*}
	\underset{T\in\mathcal{T}_h}{\text{max}}\, \text{card} (\mathcal{F}_T)\leq N_{\partial}.
\end{align*}

Next, we define the hybrid discrete spaces on the decomposition $\cT_h$. For any bounded domain $S$, let $P^k(S)$ denote the space of polynomials defined on $S$ of degree at most $k \geq 0$. Let $P^k(\cT_h)$ denote the piecewise $k$ degree polynomial functions defined on $\Omega$. The local degrees of freedom on each polygon $T \in \cT_h$ is given by
\begin{align*}
    \Ul_T^k:=\{\bvl_T=(\bv_T,(\bv_F)_{F \in \cF_T}) : \bv_T \in [P^k(T)]^d \text{ and } \bv_F \in [P^k(F)]^d, F \in \cF_T \}.
\end{align*}
The global degrees of freedom is given by combining the face values of $\Ul_T^k$ as
\begin{align*}
    \Ul_h^k := \{ \bvl_h=((\bv_T)_{T \in \cT_h},(\bv_F)_{F \in \cF_h})\}.
\end{align*}
The hybrid space with zero boundary condition is defined as
\begin{align*}
\Ul_{h,0}^k :=\{\bvl_h \in \Ul_h^k : \bv_F=0 \text{ for all } F \in \cF_h^b\}.
\end{align*}
The restriction of $\bvl_h \in \Ul_h^k$  on a polygon $T$ is denoted by $\bvl_T$. For $\bvl_h\in \Ul_h^k$, the broken polynomial function $\bv_h\in [P^k(\cT_h)]^d$ is defined as $\bv_h|_T:=\bv_T$. The local interpolation operator $\Il_T^k: [H^1(T)]^d \rightarrow \Ul_T^k$ is given by
\begin{align*}
\Il_T^k \bv :=((\pi_T^{k}\bv),(\pi_F^{k}\bv\rvert_F)_{F \in \cF_T}),
\end{align*}
where $\pi_T^{k}$ and $\pi_F^{k}$ are $L^2$ orthogonal projection onto $[P^k(T)]^d$ and $[P^k(F)]^d$ respectively. Similarly, the global interpolation operator $\Il_h^k \bv:[H^1(\Omega)]^d \rightarrow \Ul_h^k$ is defined as follows:
\begin{align*}
\Il_h^k \bv:=((\pi_T^{k}\bv)_{T \in \cT_h},(\pi_F^{k}\bv)_{F \in \cF_h}).
\end{align*}
Note that the projection operators are applied on vectors component-wise. The discrete pressure space is the usual piecewise polynomial space of degree $k$ with zero mean
\begin{align*}
P_h^k := \{p_h \in L^2(\Omega) : p_h|_T \in P^k(T)\} \cap L^2_0(\Omega).
\end{align*}
We recall some standard inequalities that will be used throughout the article.

\noindent
\textbf{Inverse Inequality:} 
\cite[Lemma~1.28]{Piet_Jero_HHO_Book_20} There exists a positive constant $C$ independent of the meshsize $h_T$ such that for any $v_h \in P^k(T)$ we have
\begin{align*}
\norm{\nabla \bv_h}_{L^2(T)}\leq Ch_T^{-1}\norm{\bv_h}_{L^2(T)}.
\end{align*}
\textbf{Trace Inequality: }~\cite[pp. 27]{Piet_Ern_12_DG_book} There exists a positive constant $C$ independent of the meshsize $h_T$ such that
\begin{equation*}
\|v\|_{L^2(\partial T)} \leq C (h_T^{-1/2} \|v\|_{L^2(T)} + h_T^{1/2} \|\nabla v\|_{L^2(T)}) \quad \forall v \in H^1(T).
\end{equation*}
In particular, for $v_h\in P^k(T)$ and $F \in \cF_T$, it holds
\begin{align*}
\|v_h\|_{L^2(F)} \leq C h_T^{-1/2} \|v_h\|_{L^2(T)} .
\end{align*}
\textbf{Discrete Poincar\'e inequality:}
\cite[Lemma~2.15]{Piet_Jero_HHO_Book_20} There exists a positive constant $C$ independent of $h_T$ such that for $\bvl_h \in \Ul_{h,0}^k$ we have
\begin{align}\label{discrete_poin}
\norm{\bv_h} \leq C \Big(\sum_{T \in \cT_h}\norm{\nabla \bv_T}_T^2+\sum_{T \in \cT_h}\sum_{F \in \cF_T}\frac{1}{h_F}\|\bv_F-\bv_T\|_F^2 \Big)^{1/2}.
\end{align}
\textbf{Approximation property of $L^2$ orthogonal projection: }~\cite[lemma 1.58]{Piet_Ern_12_DG_book} The $L^2$-projection $\pi_T^k$ satisfies the following approximation property: for any $\bv \in H^s(T)$ with $s \in \{1,2,..., k+1\}$
\begin{align}\label{approx_eqn}
\lvert \bv -\pi_T^{k}\bv \rvert_{H^m(T)} + h_T^{1/2}\lvert \bv -\pi_T^{k}\bv \rvert_{H^m(\partial T)} \leq C h_T^{s-m}\lvert \bv \rvert_{H^s(T)} \quad m \in \{0,1,...,s-1\}.
\end{align}

Let the outward unit normal component for a polygon $T \in \cT_h$ be denoted by $\bn_T$. Similarly, the outward unit normal for a face $F \in \cF_T$ is given by $\bn_{TF}$ such that $\bn_{TF}:=\bn_T\rvert_F$. Moreover, the normal component of the convection term on a face $F \in \cF_T$ is defined as $\b_{TF}:=\b\rvert_{F} \cdot\bn_{TF}$. The jump of a scalar-valued function $v$ on a face $F$ shared by two polygons $T_1$ and $T_2$ is given by
\begin{align*}
\jump{v}=v\rvert_{T_1}-v\rvert_{T_2}.
\end{align*}
The sign of $\jump{v}$ is adjusted according to the direction of the outward normal. For a vector-valued functions $\u=(u_1,\ldots, u_d)$ and  $\bv=(v_1,\ldots, v_d)$
$$ \nabla \bv:= \begin{pmatrix}
\frac{\partial v_i}{\partial x_j} 
\end{pmatrix}_{d\times d}, \hspace{0.5cm} \frac{\partial \bv}{\partial \bn}:=(\nabla \bv) \bn \;\; {\rm and}\;\; (\nabla \u,\nabla \bv):=\int_{\Omega}\sum_{i,j=1}^d\frac{\partial u_i}{\partial x_j} \frac{\partial v_i}{\partial x_j} \dx.$$
For a domain $S$, we define
$(\u,\bv)_S := \int_S \u \cdot \bv\dx$.
Let $|v|$ denote the modulus function of $v$. The positive and negative part of $v$ is defined as
\begin{align*}
v^\oplus :=\frac{1}{2}(|v|+v) \text{   and   }
v^\ominus:=\frac{1}{2}(|v|-v).
\end{align*}

\section{Discrete Oseen Problem with LPS Stabilisation}\label{Discrete_Oseen_problem}
In this section, we introduce the LPS stabilised discrete formulation for the Oseen problem \eqref{Oseen_eqn_weak} on the hybrid space $\Ul_{h,0}^k \times P_h^k$. This section is divided into three sub-sections. The first part defines some reconstruction operators that are essential to define the HHO method. In the second part, we discuss a generalised local projection stabilisation setup. The discrete LPS-HHO method is defined in the third subsection.

\subsection{Local Reconstructions} We define three reconstruction operators on the local spaces $\Ul_T^k$, see \cite{Oseen_HHO}. These are used to define the discrete HHO bilinear form in \eqref{Oseen_discrete_eqn}.

\textbf{Local velocity reconstruction:} The velocity reconstruction operator $\r_T^{k+1}:\Ul_T^k \rightarrow [P^{k+1}(T)]^d$ is defined as follows: For any $\bvl_T \in \Ul_T^k$, $\r_T^{k+1}\bvl_T$ must satisfy
\begin{equation}\label{velo_recon}
\begin{aligned}
(\nabla(\r_T^{k+1}\bvl_T),\nabla \w)_T&=(\nabla\bv_T, \nabla \w)_T+\sum_{F \in \cF_T}(\bv_F-\bv_T,\nabla \w \bn_{TF})_F \quad \forall  \w \in [P^{k+1}(T)]^d,\\
(\r_T^{k+1} \bvl_T,1)_T &=(\bv_T,1)_T.
\end{aligned}
\end{equation}
\emph{Approximation property of $\r_T^{k+1}$}:
There exists a real number $C>0$, depending on $\varrho$ but independent of $h_T$ such that, for all $\bv\in [H^{s+1}(T)]^d$ for some  $s\in\{0,1,\ldots,k+1\}$,
\begin{align*}
\|\bv-\r_T^{k+1}\Il_T^k \bv\|_T+ h_T^{1/2}\|\bv-\r_T^{k+1}\Il_T^k \bv\|_{\partial T}+ h_T\|\nabla(\bv-\r_T^{k+1}\Il_T^k \bv)\|_T \leq Ch_T^{s+1} \|\bv\|_{H^{s+1}(T)}.
\end{align*}
For $s \in \{1,2,...,k+1\}$ and $\bv \in [H^{s+1}(T)]^d$ we also have the approximation property
\begin{align}\label{recon_approximation_1}
\|\nabla(\bv-\r_T^{k+1}\Il_T^k \bv)\|_{\partial T}\leq Ch_T^{s-1/2} \|\bv\|_{H^{s+1}(T)}.
\end{align}

\textbf{Local advection reconstruction operator:} $\G_{\b,T}^k: \Ul_T^k \rightarrow [P^k(T)]^d$ is defined as follows: For any $\bvl_T \in \Ul_T^k$ 
\begin{align}\label{advection_reconstruction}
(\G_{\b,T}^k(\bvl_T),\w_T)_T=(\b \cdot \nabla \bv_T,\w_T)_T+\sum_{F \in \cF_T}(\b_{TF}(\bv_F-\bv_T),\w_T)_F \quad \forall \w_T \in [P^k(T)]^d.
\end{align}

\textbf{Local divergence reconstruction operator:} $D_T^k : \Ul_T^k \rightarrow P^k(T)$ is defined as follows: For any $\bvl_T \in \Ul_T^k$
\begin{align}\label{divergence_reconstruction}
(D_T^k(\bvl_T),q)_T=(\div \bv_T,q)_T+\sum_{F\in \cF_T}((\bv_F-\bv_T) \cdot \bn_{TF},q)_F \quad \forall q \in P^k(T).
\end{align}

\subsection{A Local Projection Setting} Let $\cM_h$ be a finite decomposition of the domain $\Omega$ into open subsets, possibly overlapping so that $\cup_{M \in \cM_h} \bar{M}=\bar{\Omega}$ and each $M \in \cM_h$ is a collection of $T\in \cT_h$. We assume that there exists a constant $C$ such that for any $M \in \cM_h$ the cardinality of the set $\{N \in \cM_h : N \cap M \neq \phi\}\leq C$. Let $h_M$ denote the diameter of the cell $M$. We also assume that for any cell $T  \in \cT_h$ inside $M \in \cM_h$, $h_M \lesssim h_T$.
Let $K_M: [L^2(M)]^d \rightarrow [P^{k-1}(M)]^d$ be a bounded linear operator defined by $K_M:=Id_M-\pi^{k-1}_M$, $Id_M$ being the identity map.

Let $\b_M$ to be a piecewise constant approximation of $\b$ on $M$ such that $\norm{\b_M}_{0,\infty,M} \leq C \norm{\b}_{0,\infty,M}$ and $\norm{\b -\b_M}_{0,\infty,M}\leq Ch_M \lvert\b\rvert_{1,\infty,M}$. For each cell $T$ contained in $M$, we define a local reconstruction $\G_{\b_M,T}^k: \Ul_T^k \rightarrow [P^k(T)]^d$ as follows: For $\bvl_T \in \Ul_T^k$
\begin{align*}
(\G_{\b_M,T}^k(\bvl_T),\w_T)_T=(\b_M \cdot \nabla \bv_T,\w_T)_T+\sum_{F \in \cF_T}(\b_{TF} (\bv_F-\bv_T),\w_T)_F \quad \forall \w_T \in [P^k(T)]^d.
\end{align*}
Define $\G_{\b_M,M}^k(\bvl_h)$ such that $\G_{\b_M,M}^k(\bvl_h)\rvert_{T} =\G_{\b_M,T}^k(\bvl_T)$ for each $T\subset M$ and $\G_{\b,M}^k(\bvl_h)$ as $\G_{\b,M}^k(\bvl_h)|_T=\G_{\b,T}^k(\bvl_T)$. In this article, we propose the following local projection stabilisation $A_{S,h} : \Ul_{h,0}^k \times \Ul_{h,0}^k \to \mathbb{R}$ defined by
\begin{align}\label{eqn_LPS}
A_{S,h}(\bvl_h,\wl_h)=\sum_{M \in \cM_h}\tau_M(K_M(\G_{\b_M,M}^k(\bvl_h)),K_M(\G_{\b_M,M}^k(\wl_h)))_M,
\end{align}
where $\norm{\b}_{0, \infty,M}\tau_M \leq C h_M$ with the constant $C$ independent of $M, F, h$ and the data of \eqref{Oseen_eqn_intro}, see \cite[eqn. (2.5)]{MR3003108}. We obtain an estimate for $\norm{\G_{\b,T}^k(\ul_h)-\G_{\b_M,T}^k(\ul_h)}$ as follows. Using the definition of the reconstructions $\G_{\b,T}^k$ and $\G_{\b_M,T}^k$ along with the approximation property of $\b_M$ we get
\begin{align*}
\norm{\G_{\b,T}^k(\ul_h)-\G_{\b_M,T}^k(\ul_h)}_{T}^2 &=((\b - \b_M) \nabla \u_T,\G_{\b,T}^k(\ul_h)-\G_{\b_M,T}^k(\ul_h))_T \nonumber\\
& \leq Ch_M^{1/2}\norm{\b}_{0, \infty, M}^{1/2}\norm{\b}_{1, \infty, M}^{1/2}\norm{\nabla \u_T}_{T}\norm{\G_{\b,T}^k(\ul_h)-\G_{\b_M,T}^k(\ul_h)}_T.
\end{align*}
This implies
\begin{align}\label{coercivity_17}
\norm{\G_{\b,T}^k(\ul_h)-\G_{\b_M,T}^k(\ul_h)}_{T}\leq Ch_M^{1/2}\norm{\b}_{0, \infty, M}^{1/2}\norm{\b}_{1, \infty, M}^{1/2}\norm{\nabla \u_T}_{T}.
\end{align}
\begin{remark}
Note that the decomposition $\cM_h$ can be taken to be the original decomposition $\cT_h$. The results proven in Sections \ref{Wellposedness} and \ref{AprioriEstimate} still hold with $\cM_h=\cT_h$. However, considering an overlapping decomposition $\cM_h$ can significantly decrease the number of degrees of freedom required for the local projection and makes the local projection stabilisation more robust with respect to the choice of stabilisation parameter $\tau_M$, see \cite{MR2670000}.
\end{remark}

\subsection{Discrete Formulation}
In this section, we introduce the discrete HHO-LPS method for the Oseen problem \eqref{Oseen_eqn_intro}. The discrete problem is defined as follows: Find $(\ul_h,p_h) \in \Ul_{h,0}^k \times P_h^k$ such that 
\begin{align}\label{Oseen_discrete_eqn}
A_h^{LP}((\ul_h,p_h),(\bvl_h,q_h))=(\f,\bv_h) \quad \forall  (\bvl_h,q_h) \in \Ul_{h,0}^k \times P_h^k,
\end{align}
where the combined bilinear form $A_h^{LP}(\cdot,\cdot)$ consists of the following parts:
\begin{align}\label{defn_bilin}
A_h^{LP}((\ul_h,p_h),(\bvl_h,q_h))&:=A_{\epsilon, h}(\ul_h,\bvl_h)+A_{\b,h}(\ul_h,\bvl_h)+A_{\rm{st}}((\ul_h,p_h),(\bvl_h,q_h))\nonumber\\
&\quad+B_h(\bvl_h,p_h)-B_h(\ul_h,q_h).
\end{align}

\medskip
\noindent Now, we define each of the bilinear forms introduced above.

\noindent\emph{The viscosity bilinear form $A_{\epsilon,h}$}: We use the local velocity reconstruction operator defined in \eqref{velo_recon} to define the viscosity term $A_{\epsilon,h}$. 
The local viscous bilinear form $A_{\epsilon,T} : \Ul_T^k \times \Ul_T^k \rightarrow \mathbb{R}$ is defined as
\begin{align*}
A_{\epsilon,T}(\wl_T,\bv_T)=\epsilon (\nabla \r_T^{k+1}(\wl_T),\nabla \r_T^{k+1}(\bvl_T))_T+S_{\epsilon,T}(\wl_T,\bvl_T), 
\end{align*}
where the local HHO stabilisation is defined as
\begin{align*}
&S_{\epsilon,T}(\wl_T,\bvl_T)\\
&=\frac{\epsilon}{h_T}\sum_{F \in \cF_T} \left(\pi_{F}^{k}(\w_F-\w_T-(r_T^{k+1}\wl_T-\pi_T^{k}r_T^{k+1}\wl_T)),\pi_{F}^{k}(\bv_F-\bv_T-(r_T^{k+1}\bvl_T-\pi_T^{k}r_T^{k+1}\bvl_T))\right)_F.
\end{align*}
The global HHO stabilisation term is given by $S_{\epsilon,h}(\wl_h,\bvl_h)=\sum_{T \in \cT_h}S_{\epsilon,T}(\wl_T,\bvl_T)$.
Summing over all $T \in \cT_h$ the global viscous bilinear form $A_{\epsilon,h}$ is given by
\begin{align*}
A_{\epsilon,h}(\wl_T,\bvl_T)=\sum_{T \in \cT_h}A_{\epsilon,T}(\wl_T,\bvl_T).
\end{align*}

\noindent\emph{The convection reaction bilinear form $A_{\b,h}$}: Define the local convection reaction bilinear form $A_{\b,T} : \Ul_T^k\times \Ul_T^k \rightarrow \mathbb{R}$ as follows
\begin{align}\label{eqn_convection}
A_{\b,T}(\wl_T,\bvl_T)&=-(\w_T,G_{\b,T}^k(\bv_T))_T+\sum_{F \in \cF_T}(\b_{TF}^{\ominus}(\w_F-\w_T),(\bv_F-\bv_T))_F \nonumber\\
&\quad+\sigma(\w_T,\bv_T)_T.
\end{align}
The global convective bilinear form $A_{\b,h} : \Ul_T^k \times \Ul_T^k \rightarrow \mathbb{R}$ is given by
\begin{align*}
A_{\b,h}(\wl_h,\bvl_h)=\sum_{T \in \cT_h}A_{\b,T}(\wl_T,\bvl_T).
\end{align*}

\noindent \emph{The velocity-pressure bilinear form $B_{h}$}: Using the definition of local divergence reconstruction in \eqref{divergence_reconstruction} the global velocity-pressure bilinear form $B_h : \Ul_h^k \times P_h^k \rightarrow \mathbb{R}$ is defined as
\begin{align}\label{eqn_velo_press}
B_h(\bvl_h,q_h)=-\sum_{T \in \cT_h}(D_{T}^k(\bvl_T),q_h)_T
\end{align}

\noindent \emph{Stabilisation terms}: The third term $A_{\rm{st}}(\cdot,\cdot)$ in \eqref{defn_bilin} consists of three stabilisation terms:
\begin{align*}
A_{\rm{st}}((\ul_h,p_h),(\bvl_h,q_h)):=A_{S,h}(\ul_h,\bvl_h)+A_{N,h}(\ul_h,\bvl_h)+B_{G,h}(p_h,q_h).
\end{align*}
The LPS stabilisation $A_{S,h}$ is defined in \eqref{eqn_LPS}.

\noindent \emph{Stabilisation for normal continuity}: Since the velocity functions in $\Ul_h^k$ do not provide normal continuity across faces, we enforce the following normal stabilisation:
\begin{align*}
A_{N,h}(\bvl_h,\wl_h):=\sum_{T \in \cT_h}\sum_{F \in \cF_T}((\bv_T-\bv_F) \cdot \bn_{TF},(\w_T-\w_F) \cdot \bn_{TF})_F.
\end{align*}

\noindent\emph{Pressure gradient stabilisation}: $B_{G,h}$ is to stabilise the pressure gradient defined as
\begin{align*}
B_{G,h}(q_h,r_h)&=\sum_{M \in \cM_h}\rho_M(K_M(\nabla_h q_h)),K_M(\nabla_h r_h)))_M,
\end{align*}
where $\rho_M \sim h_M$.


\section{Wellposedness of Discrete Formulation}\label{Wellposedness}
This section deals with the stability of the bilinear form $A_h^{LP}(\cdot,\cdot)$ as defined in \eqref{defn_bilin}. We consider the following norms and seminorms. 

\noindent\emph{Norms on $\Ul_{h,0}^k$}: For $\bvl_h \in \Ul_{h,0}^k$ define
\begin{align*}
\norm{\bvl_h}^2_{{1,h}}:=\sum_{T \in \cT_h}\Big(\norm{\nabla \bv_T}^2_T+\sum_{F \in \cF_T} \frac{1}{h_F} \norm{\bv_F-\bv_T}_F^2\Big),
\end{align*}
\begin{align*}
\norm{\bvl_h}^2_{\epsilon,h}:=A_{\epsilon,h}(\bvl_h,\bvl_h) ,\quad \norm{\bvl_h}^2_{\b}:=\sum_{T \in \cT_h}\Big(\sigma \norm{\bv_T}^2_T+\sum_{F \in \cF_T} \int_{F}\frac{\lvert \b_{TF}\rvert}{2} (\bv_F-\bv_T)\cdot (\bv_F-\bv_T)\Big).
\end{align*}
In the proof of the stability of our discrete scheme, we will use the fact that the norms $\epsilon^{1/2}\norm{\cdot}_{1,h}$ and $\norm{\cdot}_{\epsilon,h}$ are equivalent; see \cite{Piet_Jero_HHO_Book_20}.

\noindent\emph{Semi-norms and norms on $\Ul_{h,0}^k\times P_h^k$}: For $(\bvl_h,q_h) \in \Ul_{h,0}^k\times P_h^k$ define
\begin{align*}
\norm{(\bvl_h,q_h)}_{\rm st}^2:=A_{\rm st}((\bvl_h,q_h),(\bvl_h,q_h)), \quad \norm{(\bvl_h,q_h)}^2_{\rm supg}:=\sum_{M \in \cM_h}\gamma_M\norm{\G_{\b,M}^k(\bvl_h)+\nabla_h q_h}^2_M,
\end{align*}
where,
\begin{align*}
\gamma_M=h_M^2/(\epsilon+(1+\norm{\b}_{0,\infty,M})h_M+\sigma h_M^2).
\end{align*}
There exists a constant $C$ 
such that $\gamma_M \leq C\min\{\tau_M,\rho_M\}$. In our analysis, we consider the following combined norms on the space $\Ul_{h,0}^k\times P_h^k$:
\begin{align*}
\tnorm{(\bvl_h,q_h)}^2:=\norm{\bvl_h}_{\epsilon,h}^2+\norm{\bvl_h}_{\b}^2+\norm{(\bvl_h,q_h)}^2_{\rm{st}}.
\end{align*}
\begin{align}\label{LPS_norm}
\tnorm{(\bvl_h,q_h)}^2_{\rm LP}:=\tnorm{(\bvl_h,q_h)}^2+\frac{1}{1+\omega}\norm{(\bvl_h,q_h)}^2_{\rm{supg}}+(\epsilon+\sigma)\norm{q_h}^2,
\end{align}
where $\omega=\max_{M \in \cM_h} \frac{h_M^2 \norm{\b}_{1, \infty, M}}{\epsilon+\sigma h_M^2}$.

\begin{lemma}\label{lemma_coercivity} For given $(\ul_h,p_h) \in \Ul_{h,0}^k \times P_h^k$, the bilinear form defined in \eqref{defn_bilin} satisfies 
\begin{align*}
A^{LP}_h((\ul_h,p_h),(\ul_h,p_h)) = \tnorm{(\ul_h,p_h)}^2.
\end{align*}
\end{lemma}
\begin{proof}
Take the pair $(\bvl_h,q_h)=(\ul_h,p_h)$ as a test function in the definition of the bilinear form $A_h^{LP}(\cdot,\cdot)$ in \eqref{defn_bilin} to obtain
\begin{align}\label{coercivity_2}
A_h^{LP}((\ul_h,p_h),(\ul_h,p_h))&=A_{\epsilon,h}(\ul_h,\ul_h)+A_{\b,h}(\ul_h,\ul_h)+A_{\rm{st}}((\ul_h,p_h),(\ul_h,p_h)).
\end{align}
The first and third terms of the above equation read as
\begin{align}\label{eqt3}
A_{\epsilon,h}(\ul_h,\ul_h)+A_{\rm{st}}((\ul_h,p_h),(\ul_h,p_h))=\norm{\ul_h}^2_{\epsilon,h}+\norm{(\ul_h,p_h)}^2_{\rm{st}}.
\end{align}
For the second term of \eqref{coercivity_2}, we use the definition of $A_{\b,T}$ in \eqref{eqn_convection} and apply the integration by parts along with the assumption $\div \b=0$ to obtain (see \cite[eqn. (19)]{Oseen_HHO} for more details)
\begin{align}\label{eqt4}
A_{\b,h}(\ul_h,\ul_h)&=\sum_{T \in \cT_h}A_{\b,T}(\ul_T,\ul_T) \nonumber\\
&=\sum_{T \in \cT_h}\Big(-(\u_T,G_{\b,T}^k(\u_T))_T+\sum_{F \in \cF_T}(\b_{TF}^{\ominus}(\u_F-\u_T),(\u_F-\u_T))_F+\sigma(\u_T,\u_T)_T\Big)\nonumber\\
&=\sum_{T \in \cT_h}\Big(\sigma \norm{\u_T}^2_T+\sum_{F \in \cF_T}\frac{1}{2}\norm{\lvert\b_{TF}\rvert^{1/2}(\u_F-\u_T)}^2_F\Big)=\norm{\ul_h}_{\b}^2.
\end{align}
Combining the above expressions \eqref{coercivity_2}, \eqref{eqt3} and \eqref{eqt4}, we obtain the required result.
\end{proof}

\begin{lemma}\label{LPS_lemma}
For any given $(\ul_h,p_h) \in \Ul_{h,0}^k \times P_h^k$, there exists $\zl_h \in \Ul_{h,0}^k$ such that
\begin{align}\label{LPS_lemma_eqn}
A_h^{LP}((\ul_h,p_h),(\zl_h,0))&\geq \frac{1}{8}\sum_{M \in \cM_h}\gamma_M \norm{\G_{\b,M}^k(\ul_h)+\nabla_h p_h}_M^2-C_1(1+\omega)\tnorm{(\ul_h,p_h)}^2,
\end{align}
for some positive constant $C_1$ independent of $h$ and $\epsilon$.
\end{lemma}
\begin{proof}
We follow some steps of \cite[Theorem 4.1]{MR3003108} for the proof of this lemma. For given $(\ul_h,p_h) \in \Ul_{h,0}^k \times P_h^k$, let $\x_M=\G_{\b,M}^k(\ul_h)+\nabla_h p_h$ for $M \in \cM_h$. We define $\z_M=\gamma_M\pi_M^{k-1}\x_M$ on $M$ and extend to $\Omega$ by zero. Using the boundedness of $\pi_M^{k-1}$ we have
\begin{align}\label{eqt6}
\norm{\z_M}_M &\leq \gamma_M\norm{\x_M}_M .
\end{align}
Now define the global function $\zl_h=((\z_T)_{T\in\cT_h},(\z_F)_{F \in \cF_h})$, where $\z_T=\sum_{M\supset T} \z_M|_T$ and $z_F=\boldsymbol{0}$ for all $F \in \cF_h$. From the properties of the decomposition $\cM_h$ and using the inverse inequality, the trace inequality, and the equivalence $h_M \lesssim h_F$, we have the following bounds:
\begin{align*}
\norm{\z_h}^2 \leq \sum_{M \in \cM_h}\norm{\z_h}^2_M\leq C\sum_{M \in \cM_h}\norm{\z_M}^2_M,
\end{align*}
\begin{align}\label{eqt8}
\sum_{T\in \cT_h}\norm{\nabla \z_T}^2_T \leq C\sum_{M\in \cM_h}\frac{1}{h_M^2}\norm{\z_M}^2_M, \quad \sum_{T \in \cT_h}\sum_{F \in \cF_T}\frac{1}{h_{F}}\norm{\z_T}^2_F\leq C \sum_{M\in \cM_h}\frac{1}{h_M^2}\norm{\z_M}^2_M.
\end{align}
Taking $(\zl_h,0)$ as a test function in the bilinear form $A_h^{LP}$, we have
\begin{align*}
A_h^{LP}((\ul_h,p_h),(\zl_h,0))=A_{\epsilon,h}(\ul_h,\zl_h)+A_{\b,h}(\ul_h,\zl_h)+B_h(\zl_h,p_h)+A_{\rm{st}}((\ul_h,p_h),(\zl_h,0)).
\end{align*}
Applying the Cauchy-Schwarz inequality on $A_{\epsilon,h}(\ul_h,\zl_h)$ along with the equivalence of the norms $\epsilon^{1/2}\norm{\cdot}_{1,h}$ and $\norm{\cdot}_{\epsilon,h}$, we get
\begin{align}\label{eqt10}
A_{\epsilon,h}(\ul_h,\zl_h) \leq A_{\epsilon,h}(\ul_h,\ul_h)^{1/2}A_{\epsilon,h}(\zl_h,\zl_h)^{1/2} \leq CA_{\epsilon,h}(\ul_h,\ul_h)^{1/2} \epsilon^{1/2}\norm{\zl_h}_{1,h}.
\end{align}
Now, using the definition of the norm, the relations \eqref{eqt8} and the fact that $\z_F=0$, we have
\begin{align}\label{eqt11}
\epsilon \norm{\zl_h}_{1,h}^2=\sum_{T \in \cT_h}\Big(\epsilon\norm{\nabla \z_T}^2_T+\sum_{F \in \cF_T}\frac{\epsilon}{h_F}\norm{\z_F-\z_T}_F^2\Big) \leq C\sum_{M\in\cM_h}\frac{\epsilon}{h_M^2}\norm{\z_M}^2_M.
\end{align}
Using \eqref{eqt11} in \eqref{eqt10} along with \eqref{eqt6} and $\gamma_M\leq \frac{h_M^2}{\epsilon}$ we have
\begin{align}\label{coercivity_11}
A_{\epsilon,h}(\ul_h,\zl_h) \leq \frac{1}{16}\sum_{M\in\cM_h}\gamma_M\norm{\x_M}^2_M+CA_{\epsilon,h}(\ul_h,\ul_h).
\end{align}
Using the definition of $G_{b,T}^k$ in the bilinear form $A_{\b,h}(\cdot,\cdot)$ of \eqref{eqn_convection} 
and applying the integration by parts, we get (see \cite[eqn. (18)]{Oseen_HHO})
\begin{align}\label{coercivity_12}
A_{\b,h}(\ul_h,\zl_h)=\sum_{T \in \cT_h}\Big((\G_{\b,T}^k(\ul_T),\z_T)_T+\sigma(\u_T,\z_T)_T+\sum_{F \in \cF_T}(\b_{TF}^{\oplus}(\z_F-\z_T),(\u_F-\u_T))_F\Big).
\end{align}
Using the definition of the bilinear form $B_h$ in \eqref{eqn_velo_press} with integration by parts and $\z_F=0$, we have
\begin{align}\label{coercivity_13}
B_h(\zl_h,p_h)=\sum_{T \in \cT_h}\Big((\z_T,\nabla p_h)-\sum_{F \in \cF_T}(\z_F \cdot \bn_{TF},p_h)_F\Big)=\sum_{T \in \cT_h}(\z_T,\nabla_h p_h)_T.
\end{align}
Adding \eqref{coercivity_12} and \eqref{coercivity_13}, we obtain
\begin{align}\label{coercivity_14}
&A_{\b,h}(\ul_h,\zl_h)+B_h(\zl_h,p_h) \nonumber\\
&=\sum_{T \in \cT_h}\Big((G_{\b,T}^k(\ul_T)+\nabla_h p_h,\z_T)_T+\sigma(\u_T,\z_T)_T+\sum_{F \in \cF_T}(\b_{TF}^{\oplus}(\z_F-\z_T),(\u_F-\u_T))_F\Big).
\end{align}
Using the fact that $\x_M|_T=G_{\b,T}^k(\ul_T)+\nabla p_T$, the definition of $\z_M$ and relation \eqref{eqt6} the first term in the summation of \eqref{coercivity_14} becomes
\begin{align}\label{coercivity_15}
\sum_{T\in \cT_h}(G_{\b,T}^k(\ul_T)+\nabla p_T,\z_T)_T&=\sum_{M\in\cM_h}(\x_M,\z_M)_M \nonumber\\
&=\sum_{M\in\cM_h}(\x_M,\gamma_M \x_M)_M+(\x_M,\z_M-\gamma_M \x_M)_M \nonumber\\
&=\sum_{M\in\cM_h}\gamma_M\norm{\x_M}^2_M-\gamma_M(\x_M,K_M\x_M)_M. \nonumber\\
&\geq \frac{1}{2}\sum_{M\in\cM_h}\gamma_M\norm{\x_M}^2_M-C\sum_{M\in\cM_h}\gamma_M\norm{K_M(\x_M)}_M^2.
\end{align}
Now we estimate the second term in \eqref{coercivity_15}. This term can be controlled by applying the triangle inequality and then adding and subtracting the reconstruction $G_{{\b_M},M}$ with $\gamma_M \leq \min\{\tau_M, \rho_M\}$ as follows: 
\begin{align}\label{coercivity_16}
&\sum_{M\in\cM_h}\gamma_M\norm{K_M(\x_M)}_M^2 \nonumber\\
&\leq 2\sum_{M\in\cM_h}\gamma_M\norm{K_M(\G_{\b,M}^k(\ul_h))}_M^2+2\sum_{M\in\cM_h}\gamma_M\norm{K_M(\nabla_h p_h)}_M^2 \nonumber\\
&\leq 2\sum_{M\in\cM_h}\tau_M\norm{K_M(\G_{\b_M,M}^k(\ul_h))}_M^2+2\sum_{M\in\cM_h}\rho_M\norm{K_M(\nabla_h p_h)}_M^2 \nonumber\\
&\quad+\sum_{M\in\cM_h}2\gamma_M\norm{K_M(\G_{\b,M}^k(\ul_h)-\G_{\b_M,M}^k(\ul_h))}_M^2. \nonumber\\
& \leq 2 \tnorm{(\ul_h,p_h)}^2+\sum_{M\in\cM_h}2\gamma_M\norm{K_M(\G_{\b,M}^k(\ul_h)-\G_{\b_M,M}^k(\ul_h))}_M^2.
\end{align}
To control the second term of \eqref{coercivity_16}, we use the boundedness of the operator $K_M$, the inverse inequality, the equivalence of the norms $\norm{\cdot}_{\epsilon,h}$ and $\epsilon^{1/2}\norm{\cdot}_{1,h}$, the estimate \eqref{coercivity_17} and $\tau_M \norm{b}_{0,\infty,M} \lesssim h_M$, to get
\begin{align}\label{coercivity_18}
&\sum_{M \in \cM_h}\gamma_M\norm{K_M(\G_{\b,M}^k(\ul_h)-\G_{\b_M,M}^k(\ul_h))}_M^2 \nonumber\\
&\leq \sum_{M \in \cM_h} \gamma_M\norm{\G_{\b,M}^k(\ul_h)-\G_{\b_M,M}^k(\ul_h)}_M^2 \nonumber\\
&\leq C\sum_{M \in \cM_h} \tau_M h_M \norm{\b}_{0, \infty, M}\norm{\b}_{1, \infty, M}\norm{\nabla_h \u_h}_{M}^2\nonumber\\
& \leq C\sum_{M \in \cM_h}\Big(\frac{h_M^2\norm{\b}_{1, \infty, M}}{\epsilon+\sigma h_M^2}\sum_{T \subset M}(\epsilon+\sigma h_M^2)\norm{\nabla \u_T}_{T}^2\Big)\nonumber\\
&\leq C\max_{M \in \cM_h} \frac{h_M^2 \norm{\b}_{1, \infty, M}}{\epsilon+\sigma h_M^2}\Big(\norm{\ul_h}_{\epsilon,h}^2+\norm{\ul_h}_{\b}^2\Big)\leq C\omega\tnorm{(\ul_h,p_h)}^2.
\end{align}
Combining \eqref{coercivity_16} and \eqref{coercivity_18} and putting in \eqref{coercivity_15}, we obtain
\begin{align}\label{coercivity_19}
\sum_{T\in \cT_h}(G_{\b,T}^k(\ul_T)+\nabla_h p_h,\z_T)_T \geq \frac{1}{2}\sum_{M\in\cM_h}\gamma_M\norm{\x_M}^2_M-C(1+\omega)\tnorm{(\ul_h,p_h)}^2.
\end{align}
Using the Cauchy--Schwarz inequality, the relation \eqref{eqt6} and the fact that $\gamma_M \leq 1/\sigma$, the second term of \eqref{coercivity_14} gives
\begin{align}\label{coercivity_20}
\sum_{T \in \cT_h}\sigma(\u_T,\z_T)_T &\leq C\Big(\sum_{M \in \cM_h}\sigma\gamma_M^2\norm{\x_M}^2_M\Big)^{1/2}\Big(\sigma\sum_{T\in \cT_h}\norm{\u_{T}}_T^2\Big)^{1/2}\nonumber\\
& \leq \frac{1}{8}\sum_{M\in\cM_h}\gamma_M\norm{\x_M}^2_M+C\tnorm{(\ul_h,p_h)}^2.
\end{align}
Using $\z_F=0$ on each edge, $\gamma_M \leq h_M \norm{\b}_{0,\infty,M}$ and \eqref{eqt8}, the third term in \eqref{coercivity_14} gives
\begin{align}\label{coercivity_21}
&\sum_{T \in \cT_h}\sum_{F \in \cF_T}(\b_{TF}^{\oplus}(\u_F-\u_T),\z_F-\z_T)_{F}\nonumber\\
&\leq C(\sum_{T \in \cT_h}\sum_{F \in \cF_T}\norm{\b_{TF}(\u_F-\u_T)}^2_F)^{1/2}(\sum_{M\in\cM_h}\gamma_M\norm{\x_M}^2_M)^{1/2} \nonumber\\
& \leq \frac{1}{8}\sum_{M\in\cM_h}\gamma_M\norm{\x_M}^2_M+C\tnorm{(\ul_h,p_h)}^2.
\end{align}
Combining \eqref{coercivity_19}--\eqref{coercivity_21}, the expression in \eqref{coercivity_14} can be bounded as
\begin{align}\label{coercivity_22}
A_{\b,h}(\ul_h,\zl_h)+B_h(\zl_h,p_h)\geq \frac{1}{4}\sum_{M\in\cM_h}\gamma_M\norm{\x_M}^2_M-C(1+\omega)\tnorm{(\ul_h,p_h)}^2.
\end{align}
The last term remaining in $A_h^{LP}((\ul_h,p_h),(\zl_h,0))$ is
\begin{align*}
A_{\rm{st}}((\ul_h,p_h),(\zl_h,0))=A_{S,h}(\ul_h,\zl_h)+A_{N,h}(\ul_h,\zl_h).
\end{align*}
Using the definition of reconstruction $G_{\b_M,M}^k$ and \eqref{eqt8}, we have the following estimate
\begin{align*}
\norm{\G_{\b_M,M}^k(\zl_T)}_{T}^2&=(\b_M \cdot \nabla \z_T,\G_{\b_M,M}^k(\z_T))_T+\sum_{F \in \cF_T}(\b \cdot \bn_{TF}(\z_F-\z_T),\G_{\b_M,M}^k(\z_T))_F\\
&\leq C\frac{\norm{\b}_{0, \infty,M}\gamma_M}{h_T}\norm{\x_M}_M\norm{\G_{\b_M,M}^k(\zl_T)}_T.
\end{align*}
The LPS stabilisation term can now be estimated using the last inequality and boundedness of the operator $K_M$ with $\gamma_M \leq \tau_M$ and $\tau_M \norm{\b}_{0, \infty,M} \lesssim h_M$ as follows
\begin{align}\label{coercivity_23}
A_{S,h}(\ul_h,\zl_h)&=\sum_{M \in \cM_h}\tau_M\Big(K_M(\G_{\b_M,M}^k(\ul_h)),K_M(\G_{\b_M,M}^k(\zl_h))\Big)_M \nonumber\\
&\leq (\sum_{M \in \cM_h}\tau_M\norm{K_M(\G_{\b_M,M}^k(\ul_h))}^2_M)^{1/2}(\sum_{M \in \cM_h}\tau_M\norm{K_M(\G_{\b_M,M}^k(\zl_h))}^2_M)^{1/2} \nonumber\\
&\leq\Big{(}\sum_{M \in \cM_h}\tau_M\norm{K_M(\G_{\b_M,M}^k(\ul_h))}^2_M\Big{)}^{1/2} \Big{(}\sum_{M \in \cM_h}\frac{\tau_M\norm{\b}^2_{0, \infty,M}\gamma_M^2}{h_M^2}\norm{\x_M}_M^2\Big{)}^{1/2} \nonumber\\
&\leq \frac{1}{16}\sum_{M\in\cM_h}\gamma_M\norm{\x_M}^2_M+CA_{S,h}(\ul_h,\ul_h).
\end{align}
The normal stabilisation term can be handled using $\z_F=0$, $\gamma_M\leq h_M$, and \eqref{eqt8}
\begin{align}\label{coercivity_24}
A_{N,h}(\ul_h,\zl_h)&=\sum_{T \in \cT_h}\sum_{F \in \cF_T}((\u_T-\u_F) \cdot \bn_{TF},(\z_T) \cdot \bn_{TF})_F \nonumber\\
&\leq \frac{1}{16}\sum_{M\in\cM_h}\gamma_M\norm{\x_M}^2_M+CA_{N,h}(\ul_h,\ul_h).
\end{align}
Combining the inequalities \eqref{coercivity_11}, \eqref{coercivity_22}--\eqref{coercivity_24} we finally get \eqref{LPS_lemma_eqn}.
\end{proof}


\begin{lemma}\label{pressure_lemma}
For given $(\ul_h,p_h) \in \Ul_{h,0}^k \times P_h^k$, there exists $\tilde{\zl}_h \in \Ul_{h,0}^k$ such that 
\begin{align}\label{pressure_lemma_eqn}
A_h^{LP}((\ul_h,p_h),(\tilde{\zl}_h,0))\geq \frac{1}{2}\norm{p_h}^2-\tilde{C}\Big(\tnorm{(\ul_h,p_h)}^2+\norm{\u_h}^2\Big),
\end{align}
for some positive constant $\tilde{C}$ independent of $h$ and $\epsilon$.
\end{lemma}
\begin{proof}
For any fixed $p_h \in P_h^k$, take $\tilde{\zl}_h \in \Ul_{h,0}^k$ such that
\begin{align}\label{z_properties}
D_h^k(\tilde{\zl}_h)=-p_h \quad and \quad \norm{\tilde{\zl}_h}_{1,h}\leq C\norm{p_h},
\end{align}
see~\eqref{eqdivp} for details. Taking $(\tilde{\zl}_h,0)$ as a test function, we obtain
\begin{align*}
A_h^{LP}((\ul_h,p_h),(\tilde{\zl}_h,0))=A_{\epsilon,h}(\ul_h,\tilde{\zl}_h)+A_{\b,h}(\ul_h,\tilde{\zl}_h)+B_h(\tilde{\zl}_h,p_h)+A_{\rm{st}}((\ul_h,p_h),(\tilde{\zl}_h,0)).
\end{align*}
Using the equivalence of the norms $\norm{\cdot}_{\epsilon,h}$ and $\epsilon^{1/2}\norm{\cdot}_{1,h}$, the Cauchy--Schwarz inequlity, and the bound for $\tilde{\zl}_h$ in \eqref{z_properties}, we obtain an estimate for the viscous term as
\begin{align}\label{pressure_lemma_eqn_1}
A_{\epsilon,h}(\ul_h,\tilde{\zl}_h)\leq CA_{\epsilon,h}(\ul_h,\ul_h)^{1/2}\epsilon^{1/2}\norm{\tilde{\zl}_h}_{1,h} \leq \frac{1}{6}\norm{p_h}^2+ C \norm{\ul_h}^2_{\epsilon,h}.
\end{align}
Using the definition of the bilinear form $A_{\b,T}$ in \eqref{eqn_convection}, the Cauchy--Schwarz inequality, the discrete Poincar\'e inequality \eqref{discrete_poin} and the fact that $\norm{G_{\b,T}^k(\tilde{\zl}_T)}_T \leq C \norm{\tilde{\zl}_T}_{1,h}$, we get the estimate for the advection term as
\begin{align}\label{pressure_lemma_eqn_2}
A_{\b,T}(\ul_T,\tilde{\zl}_T)&=-(\u_T,G_{\b,T}^k(\tilde{\zl}_T))_T+\sum_{F \in \cF_T}(\b_{TF}^{\ominus}(\u_F-\u_T),(\tilde{\z}_F-\tilde{\z}_T))_F+\sigma(\u_T,\tilde{\zl}_T)_T\nonumber\\
& \leq C (\norm{\ul_h}_{\b}+\|u_h\|)(\norm{\tilde{\zl}_h}_{1,h}+\|\tilde{\z}_h\|) \leq \frac{1}{6}\norm{p_h}^2+ C(\norm{\ul_h}^2_{\b}+ \norm{\u_h}^2).
\end{align}
Using the definition of advection reconstruction for the LPS term and applying the trace and Cauchy--Schwarz inequalities, we obtain $\norm{G_{\b_M,T}^k(\tilde{\zl}_T)}_T \leq C \norm{\tilde{\zl}_T}_{1,h}$. This, along with the boundedness of the operator $K_M$ and \eqref{z_properties} yields
\begin{align*}
A_{S,h}(\ul_h,\tilde{\zl}_h)\leq C (A_{S,h}(\ul_h,\ul_h))^{1/2}\norm{p_h}.
\end{align*}
Applying the Cauchy-Schwarz inequality and \eqref{z_properties} on the normal stabilisation term, we have
\begin{align*}
A_{N,h}(\ul_h,\tilde{\zl}_h)\leq C (A_{N,h}(\ul_h,\ul_h))^{1/2}\norm{p_h}.
\end{align*}
Combining the last two inequalities, we obtain
\begin{align}\label{pressure_lemma_eqn_3}
A_{\rm st}((\ul_h,p_h),(\tilde{\zl}_h,0))\leq \frac{1}{6}\norm{p_h}^2+ C \norm{(\ul_h,p_h)}^2_{\rm st}.
\end{align}
The choice of $\tilde{\zl}_h$ in \eqref{z_properties} and the definition of $B_h$ provide $B_h(\tilde{\zl}_h,p_h)=\norm{p_h}^2$. Combining this with \eqref{pressure_lemma_eqn_1}--\eqref{pressure_lemma_eqn_3}, we arrive at \eqref{pressure_lemma_eqn}.
\end{proof}

\begin{theorem}\label{infsup_theorem}
There exists $\beta >0$  independent of $h$ and $\epsilon$ such that
\begin{align*}
\sup_{(\bvl_h,q_h) \in \Ul_{h,0}^k \times P_h^k}\frac{A_h^{LP}((\ul_h,p_h),(\bvl_h,q_h))}{\tnorm{(\bvl_h,q_h)}_{LP}} \geq \beta \tnorm{(\ul_h,p_h)}_{LP} \quad \forall (\ul_h,p_h) \in \Ul_{h,0}^k \times P_h^k.
\end{align*}

\end{theorem}
\begin{proof}
Using  the discrete Poincar\'e inequality \eqref{discrete_poin}, we have the following
\begin{align*}
(\epsilon+\sigma)\norm{\u_h}^2 &\leq C\Big(\sum_{T \in \cT_h}\epsilon\norm{\nabla \u_T}_T^2+\sum_{T \in \cT_h}\sum_{F \in \cF_T}\frac{\epsilon}{h_F}\int_{F}(\u_F-\u_T)\cdot(\u_F-\u_T)\Big)+\sigma \norm{\u_h}^2\\
&\leq C\Big(\norm{\ul_h}^2_{\epsilon,h}+\sigma\norm{\u_h}^2\Big).
\end{align*}
Multiplying \eqref{pressure_lemma_eqn} by $2(\epsilon +\sigma)$ and then applying the last inequality, we have
\begin{align}\label{disc_poin_epsilon}
A_h^{LP}((\ul_h,p_h),2(\epsilon+\sigma)(\tilde{\zl}_h,0))&\geq (\epsilon+\sigma)\norm{p_h}^2-\tilde{C}_1\tnorm{(\ul_h,p_h)}^2,
\end{align}
for some positive constant $\tilde{C}_1$ independent of $\epsilon$ and $\sigma$ (as both $\epsilon$ and $\sigma$ are bounded from above).

\noindent Now taking $(\bvl_h,q_h)=(\ul_h,p_h)+\frac{1}{1+8C_1+\tilde{C}_1}\Big(\frac{8}{1+\omega}(\zl_h,0)+2(\epsilon+\sigma)(\tilde{\zl}_h,0)\Big)$ with $C_1$ and $\zl_h$ as in \eqref{LPS_lemma_eqn} and $\tilde{C}_1$, $\tilde{\zl}_h$ as in \eqref{disc_poin_epsilon}, we get
\begin{align}\label{infsup_theorem_eqn_1}
A_h^{LP}((\ul_h,p_h),(\bvl_h,q_h))&\geq \frac{1}{1+8C_1+\tilde{C}_1}\tnorm{(\ul_h,p_h)}_{LP}^2.
\end{align}
Using the properties of $\zl_h$ in \eqref{eqt6}--\eqref{eqt8}, we have 
\begin{align*}
\tnorm{(\zl_h,0)}_{LP}\leq C \norm{(\ul_h,p_h)}_{{\rm supg}}
\end{align*}
In a similar manner using \eqref{z_properties}, we have 
\begin{align*}
\tnorm{(\tilde{\zl}_h,0)}_{LP} \leq C \norm{p_h}.
\end{align*}
The triangle inequality along with the last two inequalities yields
\begin{align}\label{infsup_theorem_eqn_2}
\tnorm{(\bvl_h,q_h)}_{LP} &\leq \tnorm{(\ul_h,p_h)}_{LP}+\frac{1}{1+8C+\tilde{C}_1}\Big(\frac{8}{1+\omega}\tnorm{(\zl_h,0)}_{LP}+2(\epsilon+\sigma)\tnorm{(\tilde{\zl}_h,0)}_{LP}\Big) \nonumber\\
&\leq C\tnorm{(\ul_h,p_h)}_{LP}.
\end{align}
Hence, combining \eqref{infsup_theorem_eqn_1} and \eqref{infsup_theorem_eqn_2} the theorem follows.
\end{proof}

\section{A Priori Error Estimates}\label{AprioriEstimate}
This section deals with the a priori error analysis for the discrete solution of velocity and pressure from \eqref{Oseen_discrete_eqn}. We employ the approximation results in \eqref{approx_eqn} to compute the a priori error under the $\tnorm{\cdot}_{{\rm LP}}$ norm.  We follow various steps in consistency and convergence that are based on the 3rd Strang Lemma; see \cite{Piet_Jero_HHO_Book_20,Oseen_HHO}. 
\begin{theorem}\label{apriori_estimate}
Let $(\u,p) \in \V \times Q$ be the solution to the continuous problem \eqref{Oseen_eqn_weak} and $(\ul_h,p_h) \in \Ul_{h,0}^k \times P_h^k$ be the solution to the discrete problem \eqref{Oseen_discrete_eqn}. Assume that $\u \in \V \cap [H^{k+2}(\cT_h)]^d$ such that $\Delta \u \in [L^2(\Omega)]^d$, and $\u \in H^{k+1}(\cM_h)$. Assume $p \in H^1(\Omega) \cap H^{k+1}(\cM_h)$. Then, the following result holds
\begin{align}\label{apriori_estimate_eqn}
\tnorm{(\Il_h^k \u-\ul_h,\pi_{h}^{k}p-p_h)}_{{\rm LP}} \leq C \Big(\sum_{T \in \cT_h}\epsilon h_T^{2(k+1)}\lvert \u \rvert_{k+2,T}^2+h_T^{2k+1}(\lvert p \rvert_{k+1,T}^2+\lvert \u \rvert_{k+1,T}^2)\Big)^{1/2},
\end{align}
for some positive constant C that does not depend on $h$ and $\epsilon$.
\end{theorem}
\begin{proof}
For simplicity of notation, set $(\tilde{\ul}_h,\tilde{p}_h):=(\Il_h^k \u,\pi_{h}^{k}p) \in \Ul_{h,0}^k \times P_h^k$. Then the error is $(\tilde{\ul}_h-\ul_h,\tilde{p}_h-p_h) \in \Ul_{h,0}^k \times P_h^k$. Now applying Theorem \ref{infsup_theorem}, we have 
\begin{align}\label{apriori_estimate_eqn_1}
\tnorm{(\tilde{\ul}_h-\ul_h,\tilde{p}_h-p_h)}_{{\rm LP}} \leq \frac{1}{\beta} \sup_{(\bvl_h,q_h) \in \Ul_{h,0}^k \times P_h^k} \frac{A_h^{LP}((\ul_h-\tilde{\ul}_h,p_h-\tilde{p}_h),(\bvl_h,q_h))}{\tnorm{(\bvl_h,q_h)}_{{\rm LP}}}.
\end{align}
We estimate each of the terms in the above bilinear form $A_h^{LP}$. Using the definition of discrete problem \eqref{Oseen_discrete_eqn}, we have the following consistency error:
\begin{align}\label{apriori_estimate_eqn_2}
A_h^{LP}((\ul_h-\tilde{\ul}_h,p_h-\tilde{p}_h),(\bvl_h,q_h))=(\f,\bv_h)-A_h^{LP}((\tilde{\ul}_h,\tilde{p}_h),(\bvl_h,q_h)).
\end{align}
Since $\Delta \u \in [L^2(\Omega)]^d$ and $\bvl_h \in \Ul_{h,0}^k$, the following identity holds
\begin{align*}
\sum_{T \in \cT_h}\sum_{F \in \cF_T}\epsilon (\bv_F,\nabla \u \cdot \bn_{TF})_F=0.
\end{align*}
Multiplying the first equation in \eqref{Oseen_eqn_intro} by $\bv_h$, applying the integration by parts on $(\Delta \u, \bv_h)_T$ and using the previous identity, we obtain
\begin{align}\label{apriori_estimate_eqn_3}
(\f,\bv_h)=\sum_{T \in \cT_h}\Big(\epsilon(\nabla \u,\nabla \bv_T)_T+\sum_{F \in \cF_T}\epsilon (\bv_F-\bv_T,\nabla \u \cdot \bn_{TF})_F + (\b \cdot \nabla \u,\bv_T)_T+\sigma(\u,\bv_T)_T\Big) +(\bv_h,\nabla p).
\end{align}
Now, using the expression $(\f,\bv_h)$ in \eqref{apriori_estimate_eqn_3} and the definition of $A_h^{LP}((\tilde{\ul}_h,\tilde{p}_h),(\bvl_h,q_h))$, \eqref{apriori_estimate_eqn_2} can be rewritten as
\begin{align}\label{apriori_estimate_eqn_4}
&A_h^{LP}((\ul_h-\tilde{\ul}_h,p_h-\tilde{p}_h),(\bvl_h,q_h))\nonumber\\
&=\sum_{T \in \cT_h}\Big(\epsilon(\nabla \u, \nabla \bv_T)_T+\sum_{F \in \cF_T}\epsilon (\bv_F-\bv_T,\nabla \u \cdot \bn_{TF})_F-A_{\epsilon,T}(\tilde{\ul}_T,\bvl_T)\Big)\nonumber\\
&\quad+\sum_{T \in \cT_h}\Big((\b \cdot \nabla \u,\bv_T)_T+\sigma(\u,\bv_T)_T-A_{\b,T}(\tilde{\ul}_T,\bvl_T)\Big)\nonumber\\
&\qquad+\Big((\bv_h,\nabla p)-B_{h}(\bvl_h,\tilde{p}_h)\Big)+B_h(\tilde{\ul}_h,q_h)-A_{\rm{st}}((\tilde{\ul}_h,\tilde{p}_h),(\bvl_h,q_h)) \nonumber\\
&=:E_1+E_2+E_3+E_4-E_5.
\end{align}

We estimate each of the five above terms of \eqref{apriori_estimate_eqn_4} starting with the diffusion consistency term $E_1$.  Using the definition of the reconstruction operator $\r_T^{k+1}$ in \eqref{velo_recon} with $\w=\nabla r_T^{k+1} \tilde{\ul}_T$, we have
\begin{align}\label{apriori_estimate_eqn_5}
(\nabla \r_T^{k+1} \bvl_T, \nabla \r_T^{k+1} \tilde{\ul}_T)_T=(\nabla \bv_T,\nabla \r_T^{k+1} \tilde{\ul}_T)_T+\sum_{F \in \cF_T}(\bv_F-\bv_T,\nabla \r_T^{k+1} \tilde{\ul}_T \cdot \bn_{TF})_F.
\end{align}
Using \eqref{apriori_estimate_eqn_5}, the first summation $E_1$ in \eqref{apriori_estimate_eqn_4} can be written as
\begin{align}\label{apriori_estimate_eqn_6}
&\sum_{T \in \cT_h}\Big(\epsilon(\nabla \u,\nabla \bv_T)_T+\sum_{F \in \cF_T}\epsilon (\bv_F-\bv_T,\nabla \u \cdot \bn_{TF})_F-A_{\epsilon,T}(\tilde{\ul}_T,\bvl_T)\Big) \nonumber\\
&=\sum_{T \in \cT_h}\Big(\epsilon (\nabla (\u - \r_T^{k+1}\tilde{\ul}_T),\nabla \bv_T)_T+\sum_{F \in \cF_T} \epsilon (\nabla (\u - \r_T^{k+1}\tilde{\ul}_T)\cdot \bn_{TF},\bv_F-\bv_T)_F- S_{\epsilon,T}(\tilde{\ul}_T,\bvl_T) \Big).
\end{align}
The velocity reconstruction operator satisfies $\r_T^{k+1}\Il_T^k \u=\pi_T^{1,k+1}\u$, see \cite[Definition 1.39]{Piet_Jero_HHO_Book_20}. Using this, the first term within the summation of \eqref{apriori_estimate_eqn_6} vanishes. The second term of \eqref{apriori_estimate_eqn_6} can be controlled using the approximation property of $\r_T^{k+1}$ in \eqref{recon_approximation_1} along with the Cauchy--Schwarz inequality and equivalence of the norms $\norm{\cdot}_{\epsilon,h}$ and $\norm{\cdot}_{1,h}$ as
\begin{align*}
&\sum_{F \in \cF_T} \epsilon (\nabla (\u - \r_T^{k+1}\tilde{\ul}_T)\cdot \bn_{TF},\bv_F-\bv_T)_F \leq C \epsilon^{1/2} h_T^{(k+1)}\lvert\u\rvert_{k+2,T}\norm{\bvl_h}_{\epsilon,T}.
\end{align*}
Using this and summing over all $T \in \cT_h$, we have
\begin{align}\label{apriori_estimate_eqn_8}
\sum_{T \in \cT_h}\sum_{F \in \cF_T} \epsilon (\nabla (\u - \r_T^{k+1}\tilde{\ul}_T)\cdot \bn_{TF},\bv_F-\bv_T)_F \leq C \Big(\sum_{T \in \cT_h}\epsilon h_T^{2(k+1)}\lvert\u\rvert^2_{k+2,T}\Big)^{1/2}\norm{\bvl_h}_{\epsilon,h}.
\end{align}
Since $S_{\epsilon,h}(\tilde{\ul}_h,\tilde{\ul}_h) \leq C \sum_{T \in \cT_h} \epsilon h_T^{2(k+1)}\lvert\u\rvert^2_{k+2,T}$ (see \cite{Oseen_HHO}), the third term of \eqref{apriori_estimate_eqn_6} can be controlled by using the definition of $S_{\epsilon,h}$ and the Cauchy-Schwarz inequality as
\begin{align}\label{apriori_estimate_eqn_9}
S_{\epsilon,h}(\tilde{\ul}_h,\bvl_h)\leq S_{\epsilon,h}(\tilde{\ul}_h,\tilde{\ul}_h)^{1/2} S_{\epsilon,h}(\bvl_h,\bvl_h)^{1/2}\leq C\Big(\sum_{T \in \cT_h} \epsilon h_T^{2(k+1)}\lvert\u\rvert^2_{k+2,T}\Big)^{1/2}\norm{\bvl_h}_{\epsilon,h}.
\end{align}
Combining \eqref{apriori_estimate_eqn_6}--\eqref{apriori_estimate_eqn_9}, we have
\begin{align}\label{apriori_estimate_eqn_9.1}
E_1 \leq C\Big(\sum_{T \in \cT_h} \epsilon h_T^{2(k+1)}\lvert\u\rvert^2_{k+2,T}\Big)^{1/2}\norm{\bvl_h}_{\epsilon,h}.
\end{align}
Now we estimate the term $E_2$ in \eqref{apriori_estimate_eqn_4}. Applying an integration by parts on $(\b \cdot \nabla \u,\bv_T)_T$, using the definition of $A_{\b,T}(\tilde{\ul}_T,\bvl_T)$ and the fact that $ \sum_{T\in\cT_h}\sum_{F\in\cF_T}((\b \cdot \bn_{TF})\u, \bv_F)_F=0$, we have
\begin{align}\label{apriori_estimate_eqn_10}
&\sum_{T \in \cT_h}\Big((\b \cdot \nabla \u,\bv_T)_T+\sigma(\u,\bv_T)_T-A_{\b,T}(\tilde{\ul}_T,\bvl_T)\Big) \nonumber \\
& = \sum_{T \in \cT_h}\Big(-(\u,\b \cdot \nabla \bv_T)_T+\sum_{F \in \cF_T}((\b \cdot \bn_{TF}) \u,\bv_T)_F+\sigma(\u,\bv_T)_T-A_{\b,T}(\tilde{\ul}_T,\bvl_T)\Big) \nonumber \\
& = \sum_{T \in \cT_h}-(\u-\tilde{\u}_T,\b \cdot \nabla \bv_T)_T+\sum_{T \in \cT_h}\sum_{F \in \cF_T}(\b \cdot \bn_{TF} (\u-\tilde{\u}_T),\bv_F-\bv_T)_F \nonumber\\
& \quad -\sum_{T \in \cT_h}\sum_{F \in \cF_T}\Big((\b \cdot \bn_{TF})^{\ominus}(\tilde{\u}_F-\tilde{\u}_T),(\bv_F-\bv_T)\Big)_F +\sum_{T \in \cT_h}\sigma(\u-\tilde{\u}_T,\bv_T)_T.
\end{align}
Let $\b_T$ be a $P_0$ approximation of $\b$ on $T$. Since $\pi_T^k$ is the $L^2$ projection, $(\b_T \cdot \nabla \bv_T,\u-\tilde{\u}_T)_T=0$ for all $T \in \cT_h$. We subtract this from the first term in \eqref{apriori_estimate_eqn_10} and use \eqref{approx_eqn}, inverse inequality and $\b-\b_T\sim h_T$ to get
\begin{align}\label{apriori_estimate_eqn_11}
\sum_{T \in \cT_h}(\u-\tilde{\u}_T,(\b-\b_T) \cdot \nabla \bv_T)_T &\leq  C\Big(\sum_{T \in \cT_h}  h_T^{2k+2}\lvert\u\rvert^2_{k+1,T}\Big)^{1/2}\Big(\sum_{T \in \cT_h}\norm{\bv_T}^2_T\Big)^{1/2} \nonumber\\
&\leq C \Big(\sum_{T \in \cT_h}  h_T^{2k+2}\lvert\u\rvert^2_{k+1,T}\Big)^{1/2}\norm{\bvl_h}_{\b}.
\end{align}
Using the approximation estimates in \eqref{approx_eqn}, the Cauchy--Schwarz inequality, and the trace inequality, the second term in \eqref{apriori_estimate_eqn_10} gives
\begin{align}\label{apriori_estimate_eqn_12}
\sum_{T \in \cT_h}\sum_{F \in \cF_T}(\b \cdot \bn_{TF} (\u-\tilde{\u}_T),\bv_F-\bv_T)_F & \leq C\Big(\sum_{T \in \cT_h}  h_T^{2k+1}\lvert\u\rvert^2_{k+1,T}\Big)^{1/2}\Big(\sum_{T \in \cT_h}\sum_{F \in \cF_T}\frac{1}{2}\norm{\b \cdot \bn_{TF}(\bv_F-\bv_T)}_F^2 \Big)^{1/2} \nonumber\\
&\leq C\Big(\sum_{T \in \cT_h}  h_T^{2k+1}\lvert\u\rvert^2_{k+1,T}\Big)^{1/2}\norm{\bvl_h}_{\b}.
\end{align}
The third term of \eqref{apriori_estimate_eqn_10} is the upwind stabilisation term which can be controlled using the estimate
\begin{align*}
\| \tilde{\u}_F-\tilde{\u}_T\|_{F}=\|\pi_F^k\u-\pi_T^k\u\|_F=\|\pi_F^k(\u-\pi_T^k\u)\|_F\leq \|\u-\pi_T^k\u\|_F
\end{align*}
and steps similar to \eqref{apriori_estimate_eqn_12} as
\begin{align}\label{apriori_estimate_eqn_13}
\sum_{T \in \cT_h}\sum_{F \in \cF_T}\Big((\b \cdot \bn_{TF})^{\ominus}(\tilde{\u}_F-\tilde{\u}_T),(\bv_F-\bv_T)\Big)_F\leq C\Big(\sum_{T \in \cT_h}  h_T^{2k+1}\lvert\u\rvert^2_{k+1,T}\Big)^{1/2}\norm{\bvl_h}_{\b}.
\end{align}
The last term in \eqref{apriori_estimate_eqn_10} is the reaction term which can be simply bounded as
\begin{align}\label{apriori_estimate_eqn_14}
\sigma(\u-\tilde{\u}_T,\bv_T)_T \leq C \Big(\sum_{T \in \cT_h}  h_T^{2(k+1)}\lvert\u\rvert^2_{k+1,T}\Big)^{1/2}\norm{\bvl_h}_{\b}.
\end{align}
Combining the estimates \eqref{apriori_estimate_eqn_11}--\eqref{apriori_estimate_eqn_14}, we get
\begin{align}\label{apriori_estimate_eqn_14.1}
E_2 \leq C \Big(\sum_{T \in \cT_h}  h_T^{2k+1}\lvert\u\rvert^2_{k+1,T}\Big)^{1/2}\norm{\bvl_h}_{\b}.
\end{align}
The third term of the consistency error in \eqref{apriori_estimate_eqn_4} is $E_3=(\bv_h,\nabla p)-B_{h}(\bvl_h,\tilde{p}_h)$. Using the definition of the bilinear form $B_h$ and using $(\nabla \cdot \bv_T, \pi_T^{k}p)_T=(\nabla \cdot \bv_T, p)_T$, we arrive at
\begin{align*}
B_h(\bvl_h,\tilde{p}_h)&=\sum_{T \in \cT_h}\Big(-(\nabla \cdot \bv_T, \tilde{p}_T)_T-\sum_{F \in \cF_T}((\bv_F-\bv_T)\cdot \bn_{TF},\tilde{p}_T)_F\Big) \nonumber\\
&=\sum_{T \in \cT_h}\Big(-(\nabla \cdot \bv_T, p)_T-\sum_{F \in \cF_T}((\bv_F-\bv_T)\cdot \bn_{TF},\tilde{p}_T)_F\Big).
\end{align*}
Since $p \in H^1(\Omega)$ and $\bv_F=0$ on $\partial \Omega$, we have $\sum_{T \in \cT_h}\sum_{F \in \cF_T}(\bv_F \cdot \bn_{TF},p)_F=0$. Using this, along with an integration by parts on $(\bv_T,\nabla p)_T$, we obtain
\begin{align*}
(\bv_h,\nabla p)=\sum_{T \in \cT_h}\Big(-(\nabla \cdot \bv_T,p)_T-\sum_{F \in \cF_T}((\bv_F-\bv_T)\cdot \bn_{TF},p)_F\Big).
\end{align*}
Combining the last two inequalities, $E_3$ becomes
\begin{align}\label{apriori_estimate_eqn_15}
E_3=(\bv_h,\nabla p)-B_{h}(\bvl_h,\tilde{p}_h)&=\sum_{T \in \cT_h}\sum_{F \in \cF_T}((\bv_F-\bv_T)\cdot \bn_{TF},\tilde{p}-p)_F \nonumber\\
&\leq C\Big(\sum_{T \in \cT_h}h_T^{2k+1}\lvert p \rvert^2_{k+1,T} \Big)^{1/2}\Big(\sum_{T \in \cT_h}\sum_{F \in \cF_T}\norm{(\bv_F-\bv_T)\cdot \bn_{TF}}^2_F\Big)^{1/2}\nonumber\\
&\leq C\Big(\sum_{T \in \cT_h}h_T^{2k+1}\lvert p \rvert^2_{k+1,T} \Big)^{1/2}\Big(A_{N,h}(\bvl_h,\bvl_h)\Big)^{1/2}.
\end{align}
The fourth term $E_4$ in \eqref{apriori_estimate_eqn_4} is $B_h(\tilde{\ul}_h,q_h)$. This term can be proved to be zero using the fact that $D_T^k\Il_T^k \u=\pi_T^{k}(\div \u)$ and $\div \u=0$
\begin{align}\label{apriori_estimate_eqn_16}
E_4=B_h(\tilde{\ul}_h,q_h)&=-\sum_{T \in \cT_h}(D_T^{k}(\Il_T^{k}\u),q_h)_T=-\sum_{T \in \cT_h}(\pi_T^{k}(\div \u),q_h)_T=0.
\end{align}
The last term of \eqref{apriori_estimate_eqn_4} is the stabilisation term $E_5=A_{\rm{st}}(\cdot,\cdot)$ which has the following three components
\begin{align}\label{apriori_estimate_eqn_17}
A_{\rm{st}}((\tilde{\ul}_h,\tilde{p}_h),(\bvl_h,q_h))=A_{S,h}(\tilde{\ul}_h,\bvl_h)+A_{N,h}(\tilde{\ul}_h,\bvl_h)+B_{G,h}(\tilde{p}_h,q_h).
\end{align}
The first term of \eqref{apriori_estimate_eqn_17} is the LPS stabilisation defined in \eqref{eqn_LPS}. Using the orthogonality of $\pi_M^{k-1}$, we have
$$A_{S,h}(\tilde{\ul}_h,\bvl_h)=\tau_M(G_{\b_M,M}(\Il_T^k(\u)),K_M(\G_{\b_M,M}^k(\bvl_h)))_M.$$
We use some intermediate steps to estimate the above term. Applying the integration by parts (twice) and using the projection property of $\pi_T^k$, we obtain for $T \subset M$
\begin{align}
&(\b_M\cdot \nabla \pi_T^k \u,K_M(\G_{\b_M,M}^k(\bvl_h)))_T\nonumber\\
&=-(\pi_T^k \u,\b_M \cdot \nabla K_M(\G_{\b_M,M}^k(\bvl_h)))_T+\sum_{F \in \cF_T}((\b_M \cdot \bn_{TF}) \pi_T^k \u,K_M(\G_{\b_M,M}^k(\bvl_h)))_F\notag\\
&=-(\u,\b_M \cdot \nabla K_M(\G_{\b_M,M}^k(\bvl_h)))_T+\sum_{F \in \cF_T}((\b_M \cdot \bn_{TF}) \pi_T^k \u,K_M(\G_{\b_M,M}^k(\bvl_h)))_F\notag\\
&=(\b_M \cdot \nabla \u,K_M(\G_{\b_M,M}^k(\bvl_h)))_T+\sum_{F \in \cF_T}((\b_M \cdot \bn_{TF}) (\pi_T^k \u-\u),K_M(\G_{\b_M,M}^k(\bvl_h)))_F.\label{Gb_IBP}
\end{align}
Using the definition $G_{\b_M,T}^k(\Il_T^k(\u))$ and the above equality \eqref{Gb_IBP}, we obtain
\begin{align}\label{apriori_estimate_eqn_18}
&(G_{\b_M,T}(\Il_T^k(\u)),K_M(\G_{\b_M,M}^k(\bvl_h)))_T\nonumber\\
&=(\b_M\cdot \nabla \pi_T^k \u,K_M(\G_{\b_M,M}^k(\bvl_h)))_T+\sum_{F \in \cF_T}((\b \cdot \bn_{TF})(\pi_F^k \u-\pi_T^k \u),K_M(\G_{\b_M,M}^k(\bvl_h)))_F \nonumber\\
&=(\b_M \cdot \nabla \u,K_M(\G_{\b_M,M}^k(\bvl_h)))_T+\sum_{F \in \cF_T}((\b \cdot \bn_{TF})(\pi_F^k \u-\pi_T^k \u),K_M(\G_{\b_M,M}^k(\bvl_h)))_F \nonumber\\
&+\sum_{F \in \cF_T}((\b_M \cdot \bn_{TF})(\pi_T^k \u-\u),K_M(\G_{\b_M,M}^k(\bvl_h)))_F.
\end{align}
Summing the last equation over all $T \subset M$ and applying $(\pi_M^{k-1}(\b_M \cdot \nabla \u),K_M(\G_{\b_M,M}^k(\bvl_h)))_M=0$ along with the approximation properties of $\pi_T^{k}$, $\pi_M^{k-1}$ and $\pi_F^{k}$, we get
\begin{align*}
&\tau_M(K_M(\G_{\b_M,M}^k(\tilde{\ul}_h)),K_M(\G_{\b_M,M}^k(\bvl_h)))_M \nonumber\\
&=\tau_M  \Big((\b_M \cdot \nabla \u,K_M(\G_{\b_M,M}^k(\bvl_h)))_M+\sum_{T \subset M}\sum_{F \in \cF_T}((\b \cdot \bn_{TF})(\pi_F^k \u-\pi_T^k \u),K_M(\G_{\b_M,M}^k(\bvl_h)))_F \nonumber\\
& \qquad+\sum_{T \subset M}\sum_{F \in \cF_T}((\b_M \cdot \bn_{TF})(\pi_T^k \u-\u),K_M(\G_{\b_M,M}^k(\bvl_h)))_F\Big) \nonumber\\
&\leq C h_M^{k+1/2} \lvert \u \rvert_{k+1,M}\tau_M^{1/2}\norm{K_M(\G_{\b_M,M}^k(\bvl_h))}_M.
\end{align*}
In particular, $\tau_M^{1/2}\|K_M(\G_{\b_M,M}^k(\tilde{\ul}_h)\|_M\leq C h_M^{k+1/2} \lvert \u \rvert_{k+1,M}$.
Using this with $h_M\lesssim h_T$, the LPS stabilisation term can be bounded as
\begin{align}\label{apriori_estimate_eqn_19}
A_{S,h}(\tilde{\ul}_h,\bvl_h)\leq A_{S,h}(\tilde{\ul}_h,\tilde{\ul}_h)^{1/2} A_{S,h}(\bvl_h,\bvl_h)^{1/2}\leq C\Big(\sum_{T \in \cT_h}  h_T^{2k+1}\lvert\u\rvert^2_{k+1,T}\Big)^{1/2} A_{S,h}(\bvl_h,\bvl_h)^{1/2}.
\end{align}
The normal jump stabilisation term $A_{N,h}$ can be controlled using the approximation property of $\pi_T^{k}$ in \eqref{approx_eqn} and the boundedness of $\pi_F^{k}$ as follows
\begin{align}\label{apriori_estimate_eqn_20}
A_{N,h}(\tilde{\ul}_h,\bvl_h)&=\sum_{T \in \cT_h}\sum_{F\in \cF_T}\int_{F}((\pi_F^k\u-\pi_T^k \u)\cdot \bn_{TF},(\bv_F-\bv_T)\cdot \bn_{TF})_F \nonumber\\
&=\sum_{T \in \cT_h}\sum_{F\in \cF_T}\int_{F}(\pi_F^k(\u-\pi_T^k \u)\cdot \bn_{TF},(\bv_F-\bv_T)\cdot \bn_{TF})_F \nonumber\\
& \leq C\Big(\sum_{T \in \cT_h}h_T^{2k+1}\lvert \u \rvert^2_{k+1}\Big)^{1/2}A_{N,h}(\bvl_h,\bvl_h)^{1/2}.
\end{align}
The approximation property of $K_M$ gives $\norm{K_M(\nabla p)}_M \leq C h_M^k \norm{p}_{k+1,M}$. Using this and the boundedness of the operator $K_M$ along with the approximation properties of $\pi_T^{k}$, we get
\begin{align}\label{apriori_estimate_eqn_21}
B_{G,h}(\tilde{p}_h,q_h)&\leq B_{G,h}(\tilde{p}_h,\tilde{p}_h)^{1/2} B_{G,h}(q_h,q_h)^{1/2} \nonumber\\
& \leq \Big(\rho_M^{1/2}\norm{K_M(\nabla (p-\tilde{p}_h))}_M+\rho_M^{1/2}\norm{K_M(\nabla p)}_M\Big)B_{G,h}(q_h,q_h)^{1/2} \nonumber\\
& \leq C\Big(\sum_{T \in \cT_h}h_T^{2k+1}\lvert p \rvert_{k+1}^2\Big)^{1/2}B_{G,h}(q_h,q_h)^{1/2}.
\end{align}
Combining \eqref{apriori_estimate_eqn_19}--\eqref{apriori_estimate_eqn_21}, the stabilisation term in \eqref{apriori_estimate_eqn_4} can be bounded as
\begin{align}\label{apriori_estimate_eqn_23}
A_{\rm{st}}((\tilde{\ul}_h,\tilde{p}_h),(\bvl_h,q_h))\leq C\Big(\sum_{T \in \cT_h}h_T^{2k+1}\lvert \u \rvert^2_{k+1}+h_T^{2k+1}\lvert p \rvert_{k+1}^2\Big)^{1/2}A_{\rm{st}}((\bvl_h,q_h),(\bvl_h,q_h))^{1/2}.
\end{align}
Finally, combining \eqref{apriori_estimate_eqn_1}, \eqref{apriori_estimate_eqn_9.1},\eqref{apriori_estimate_eqn_14.1}--\eqref{apriori_estimate_eqn_16} and \eqref{apriori_estimate_eqn_23}, we obtain \eqref{apriori_estimate_eqn}.
\end{proof}

\begin{remark}
The analysis performed in \eqref{apriori_estimate_eqn_16} shows that the normal jump stabilisation term $A_{N,h}(\cdot,\cdot)$ is essential for error analysis. Comparing the analysis of this term with \cite[pg-1332 ]{Oseen_HHO} shows that there is a presence of $\epsilon$ in the denominator. In our analysis, we have bypassed this by taking the normal jump stabilisation term. 
\end{remark}
\begin{remark}
Note that if $\cM_h=\cT_h$ then the regularity assumption on the velocity and pressure spaces can be taken to be $[H^{k+1}(\cT_h)]^d$ and $H^{k+1}(\cT_h)$ respectively. Moreover, for $\cM_h=\cT_h$ the pressure gradient stabilisation term $B_{G,h}(p_h,q_h)=0$ since $K_T(\nabla p_h)=0$.
\end{remark}

\section{Numerical Results}\label{numerics}
In this section, we perform some numerical experiments for the HHO approximation of the Oseen problem~\eqref{Oseen_eqn_sep} to validate the a priori results obtained in Theorem~\ref{apriori_estimate}. The experiments are performed for the case $\cM_h=\cT_h$. Let the error be denoted by $(\underline{\boldsymbol{e}}_h^{\u}, e_h^{p}):=(\Il_h^k \u-\ul_h,\pi_{h}^{k}p-p_h)$. In the following numerical experiment, we compute the error w.r.t. the LPS norm $\tnorm{(\underline{\boldsymbol{e}}_h^{\u}, e_h^{p})}_{\rm LP}$ as defined in \eqref{LPS_norm}. We compute the empirical rate of convergence using the formula
\begin{align*}
	\texttt{rate}(\ell):=\log \big(\tnormn{(\underline{\boldsymbol{e}}_{h_{\ell}}^{\u}, e_{h_{\ell}}^{p})}_{\rm LP}/\tnormn{(\underline{\boldsymbol{e}}_{h_{\ell-1}}^{\u}, e_{h_{\ell-1}}^{p})}_{\rm LP} \big)/\log \big(h_\ell/h_{\ell-1} \big)\text{ for } \ell=1,2,3,\ldots,
\end{align*}
where $\tnormn{(\underline{\boldsymbol{e}}_{h_{\ell}}^{\u}, e_{h_{\ell}}^{p})}_{\rm LP}$ and $\tnormn{(\underline{\boldsymbol{e}}_{h_{\ell-1}}^{\u}, e_{h_{\ell-1}}^{p})}_{\rm LP}$ are the errors associated to the two consecutive meshsizes $h_\ell$ and $h_{\ell-1}$, respectively. We adopted some of the basic implementation methodologies for the HHO methods from \cite{Piet_Jero_HHO_Book_20,Cicu_Piet_Ern_18_Comp_HHO, Quasilin_HHO_GM_TG_TP}.

\begin{example}\label{test:smooth} 
Consider the Oseen problem \eqref{Oseen_eqn_intro} with the homogeneous Dirichlet boundary condition in the square domain $\Omega=(0, 1)^2$. In this experiment, we consider the convection term $\b=(1,1)$ and the reaction term $\sigma=1$.  We take the force function to be $\f:=-\epsilon \Delta \u + (\b \cdot \nabla) \u + \sigma \u  + \nabla p$, where the exact solution for velocity and pressure are given by
\begin{align*}
u(x,y)=\begin{bmatrix}
2x^2y(2y-1)(x - 1)^2(y - 1)\\-2xy^2(2x - 1)(x - 1)(y - 1)^2 \end{bmatrix},
\end{align*}
\begin{align*}
p(x,y)=2\cos(x)\sin(y)-2\sin(1)(1-\cos(1)).
\end{align*}
The numerical tests are performed on the triangular, Cartesian, and hexagonal mesh families; see Figure~\ref{fig:Triag_Cart_Hexa_Meshes} \cite{Quasilin_HHO_GM_TG_TP}. We consider the triangular and Cartesian mesh families of \cite{Her_Hub_FVCA_mesh_08}, and the hexagonal mesh family of \cite{Piet_Lem_HexaMesh_15}. In the left part of Figure~\ref{fig:Smooth_Conv_His_Triag_Cart_Hexa} ((a), (c) and (e)), we plot the error $\tnorm{(\underline{\boldsymbol{e}}_h^{\u}, e_h^{p})}_{\rm LP}$ of Theorem~\ref{apriori_estimate} as a function of meshsize $h$ for polynomial degrees $k=0,1,2,3$ with $\epsilon=10^{-8}$. We observe that the convergence rates for the errors are approximately $h^{k+1/2}$ for $\epsilon=10^{-8}$. In the right part of Figure~\ref{fig:Smooth_Conv_His_Triag_Cart_Hexa} ((b), (d) and (f)), we plot the error $\|(\boldsymbol{e}_h^{\u}, e_h^{p})\|_{\rm supg}$ in the supg norm of Theorem~\ref{apriori_estimate} as a function of meshsize $h$ for polynomial degrees $k=0,1,2,3$ with $\epsilon=10^{-8}$. The convergence histories show the same rate of convergence for both norms $\tnorm{\cdot}_{\rm LP}$ and $\|\cdot\|_{\rm supg}$. This validates the theoretical results obtained in Theorem~\ref{apriori_estimate}.

\begin{figure}
	\begin{center}
		\subfloat[]{\includegraphics[width=0.33\textwidth]{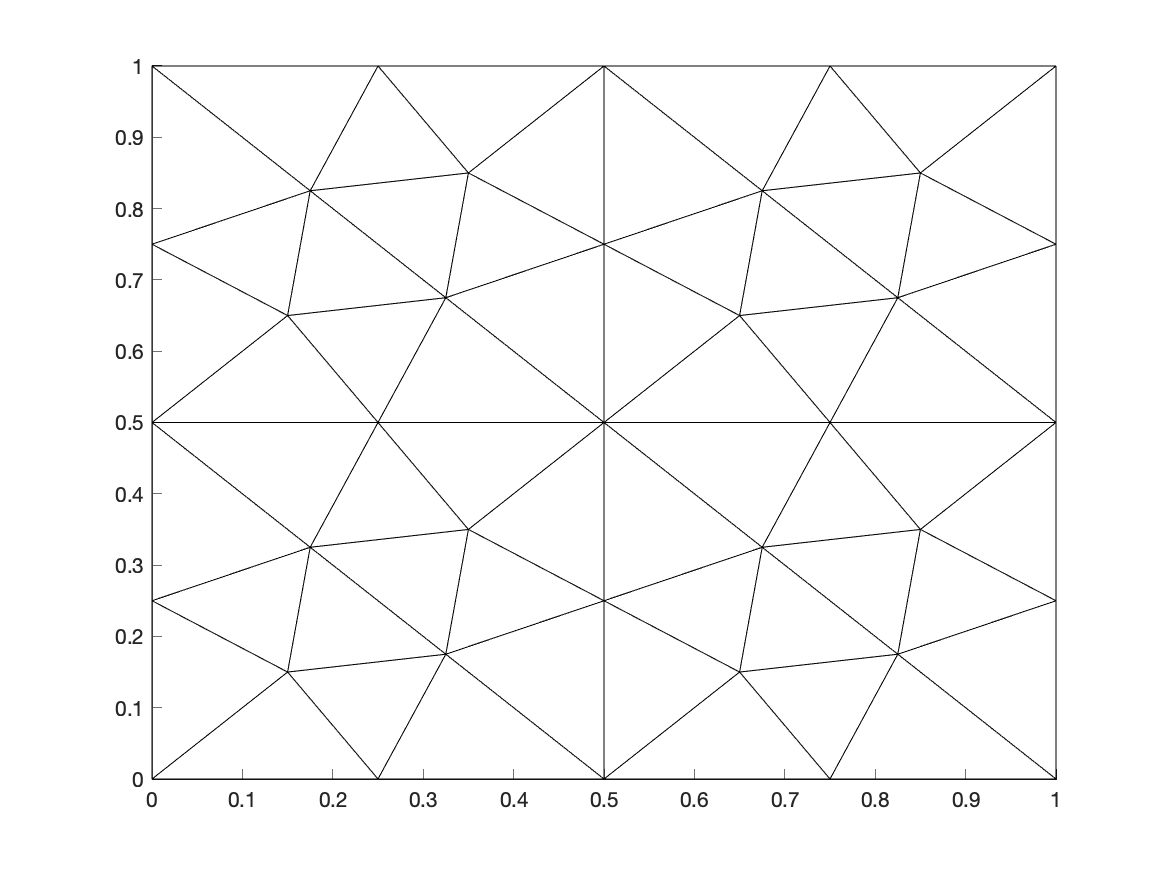}}
		\subfloat[]{\includegraphics[width=0.33\textwidth]{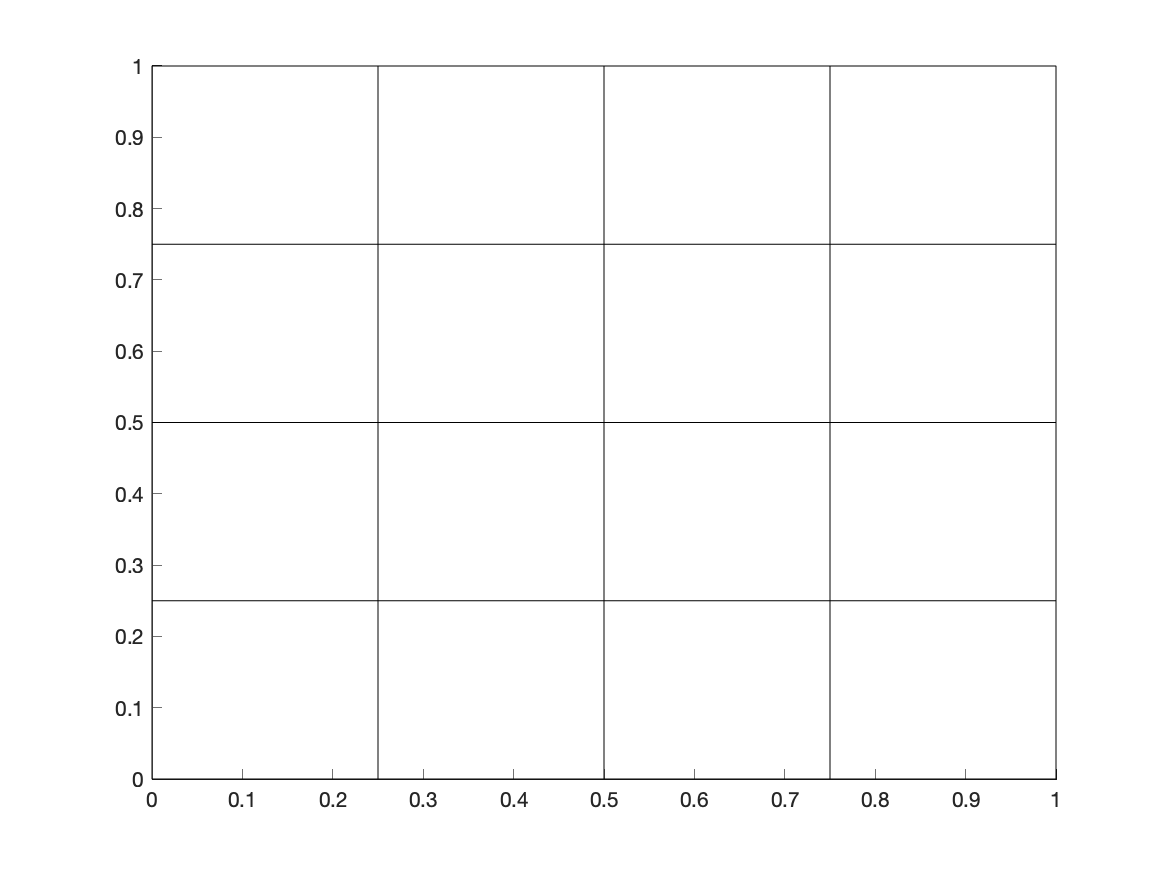}}
		\subfloat[]{\includegraphics[width=0.33\textwidth]{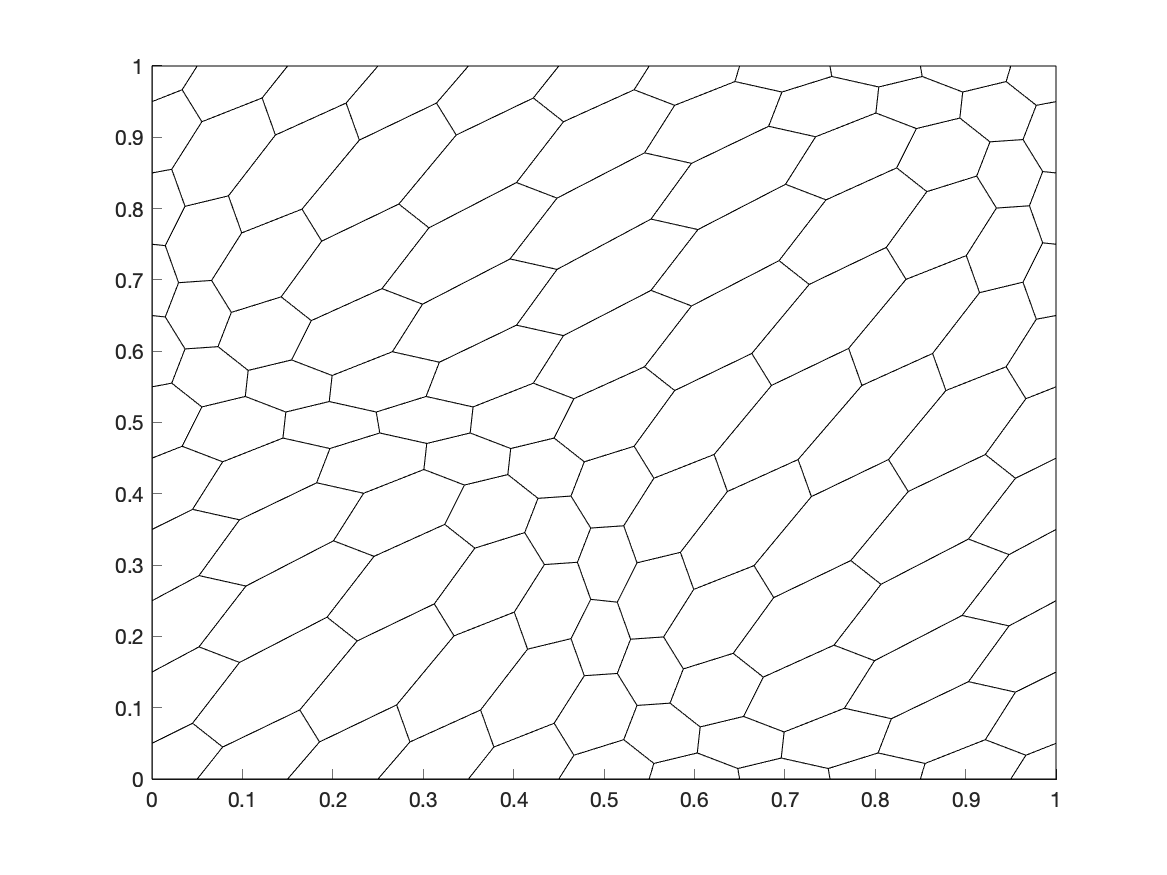}}
		\caption{(a) Triangular, (b)  Cartesian and (c) hexagonal initial meshes.}
		\label{fig:Triag_Cart_Hexa_Meshes}
	\end{center}
\end{figure}

\begin{figure}
	\begin{center}
		\subfloat[]{\includegraphics[height=0.4\textwidth,width=0.5\textwidth]{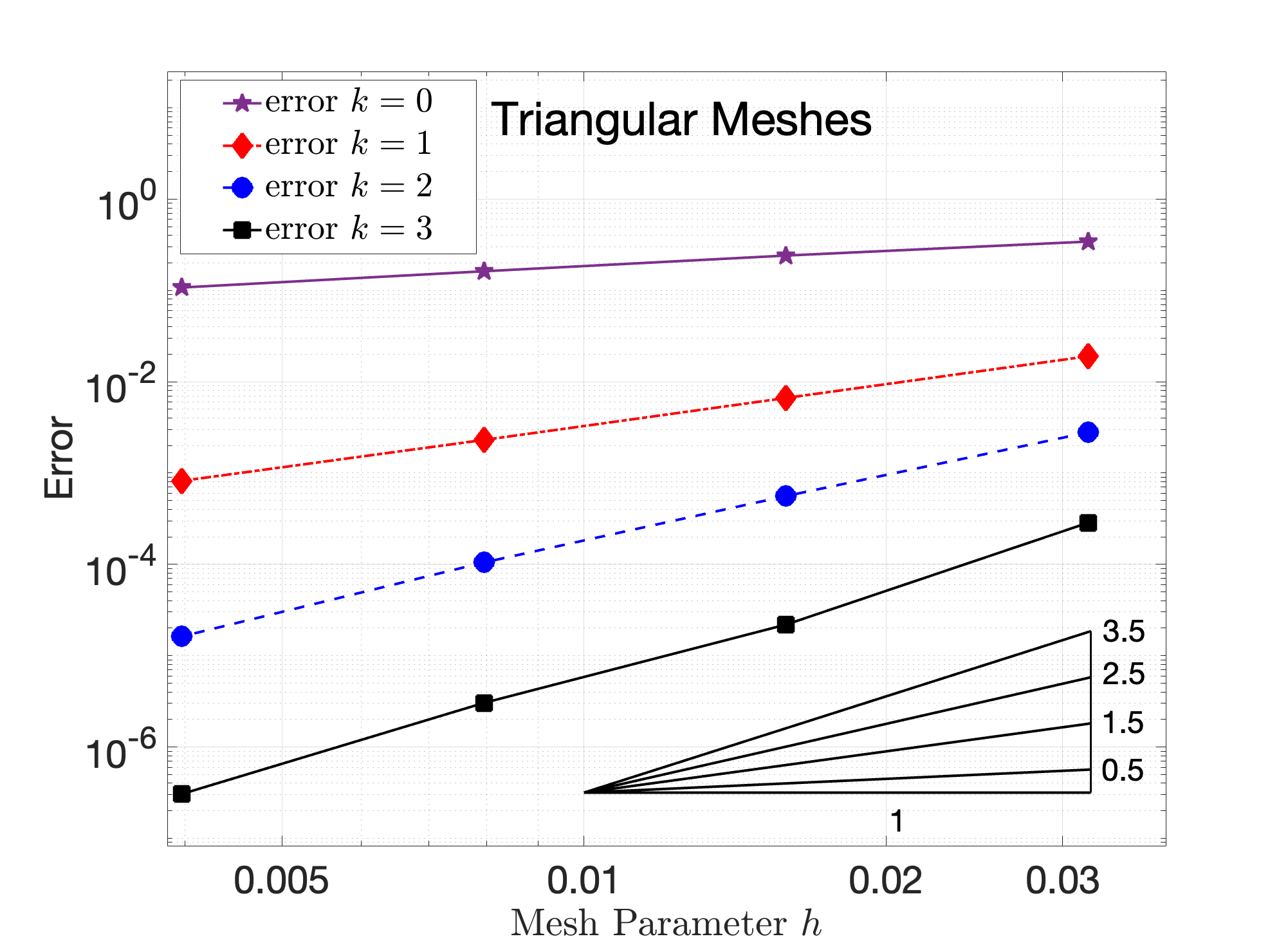}}
		\subfloat[]{\includegraphics[height=0.4\textwidth,width=0.5\textwidth]{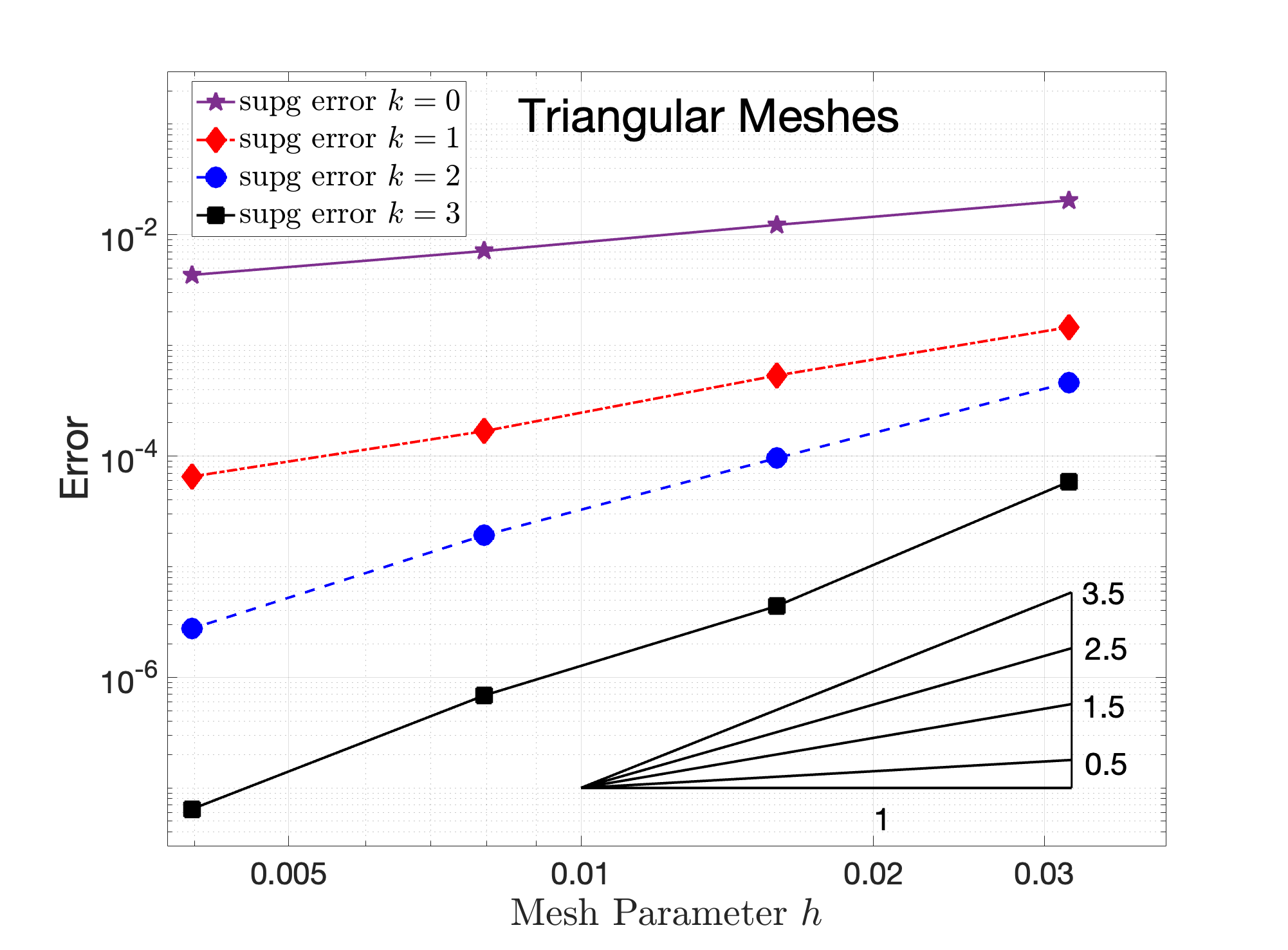}}\\
		\subfloat[]{\includegraphics[height=0.4\textwidth,width=0.5\textwidth]{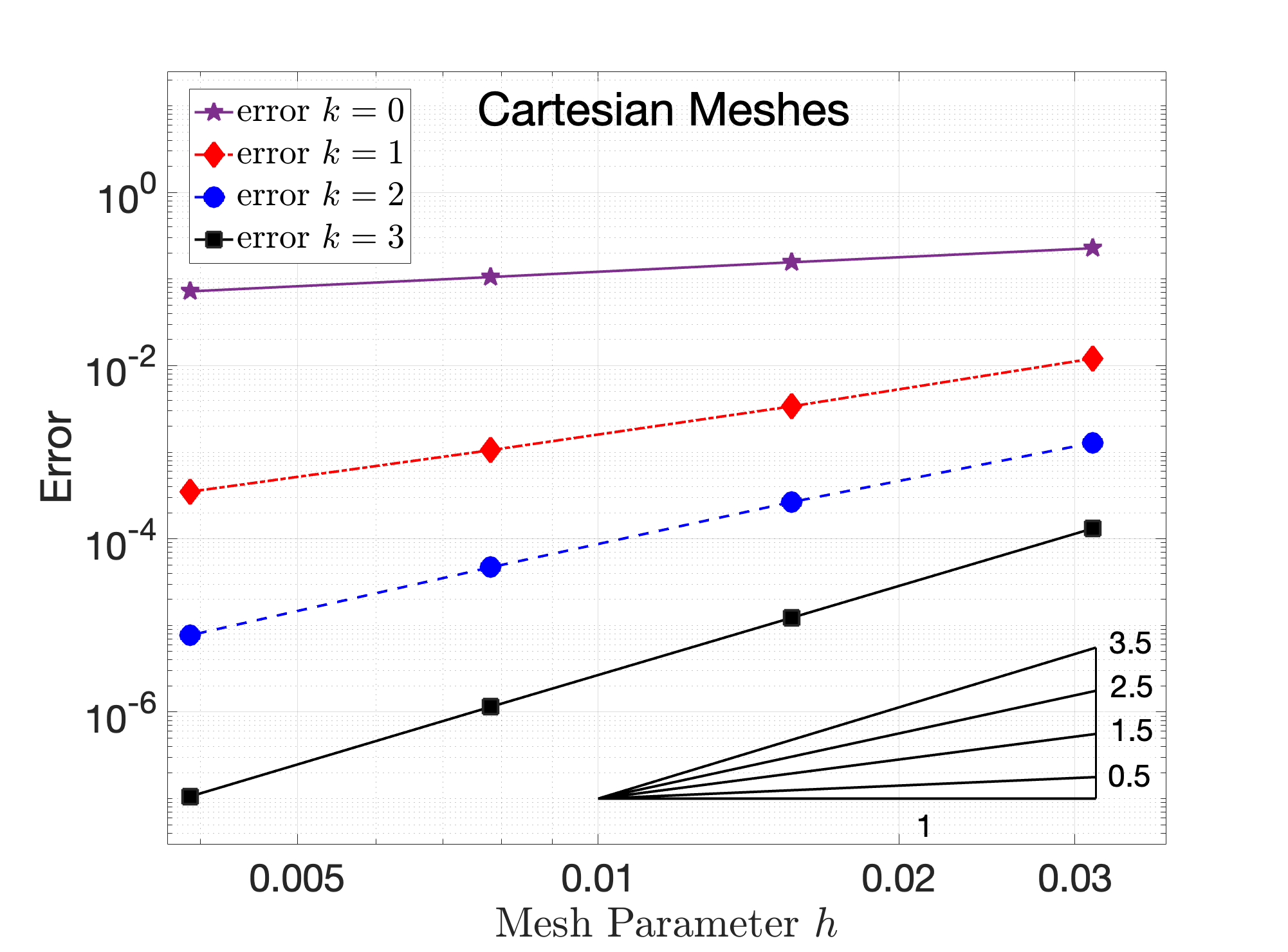}}
		\subfloat[]{\includegraphics[height=0.4\textwidth,width=0.5\textwidth]{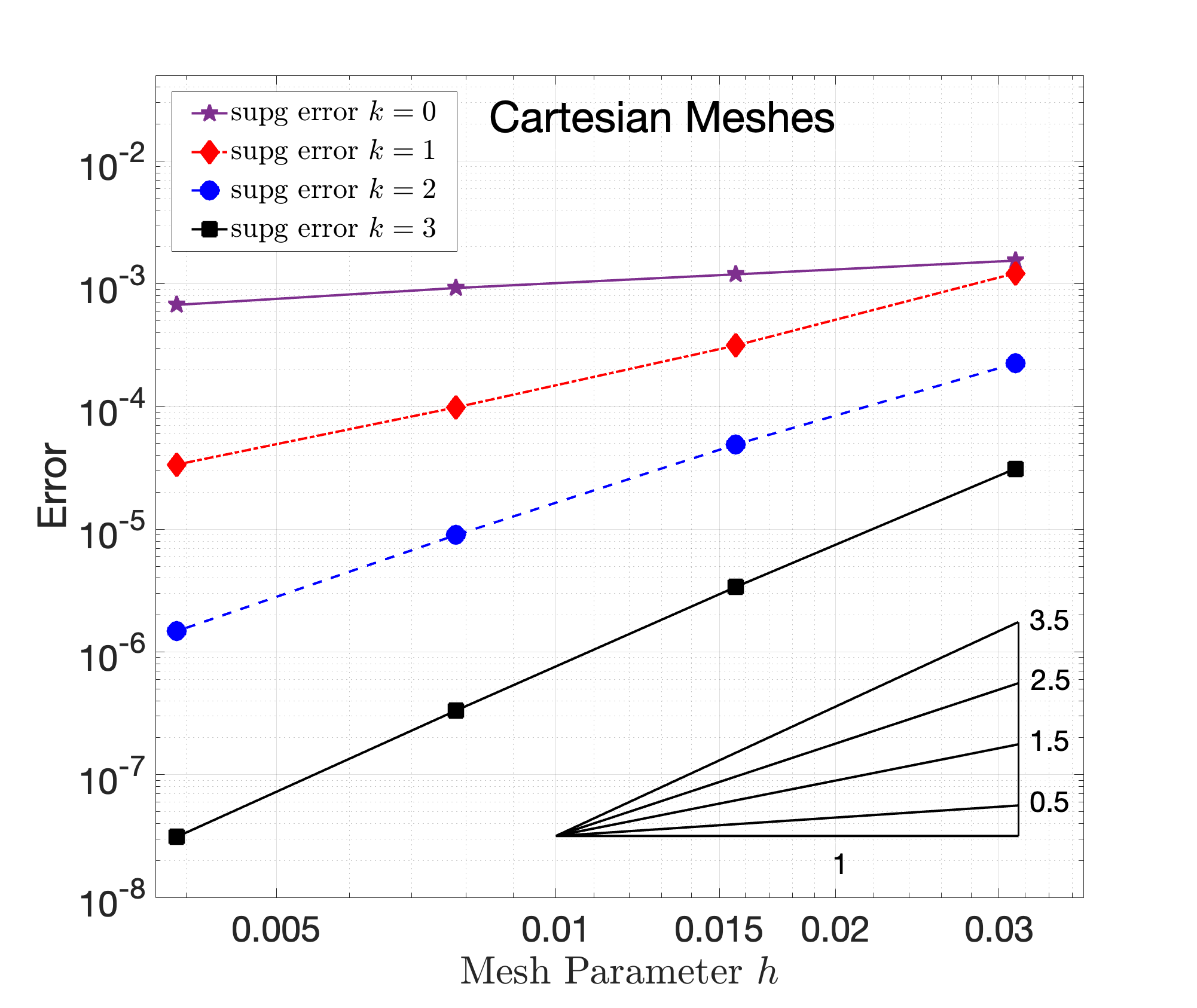}}\\
		\subfloat[]{\includegraphics[height=0.4\textwidth,width=0.5\textwidth]{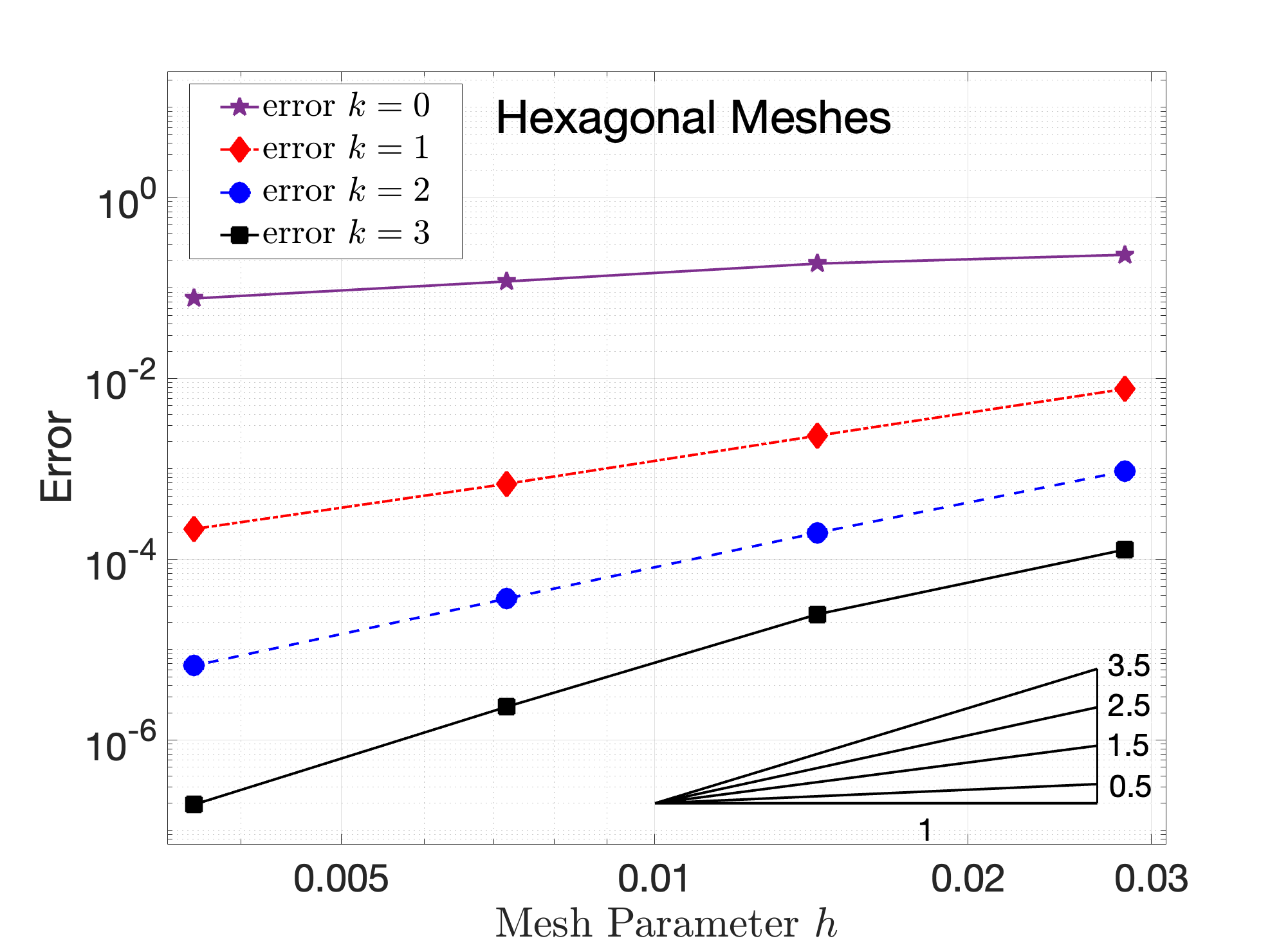}}
  		\subfloat[]{\includegraphics[height=0.4\textwidth,width=0.5\textwidth]{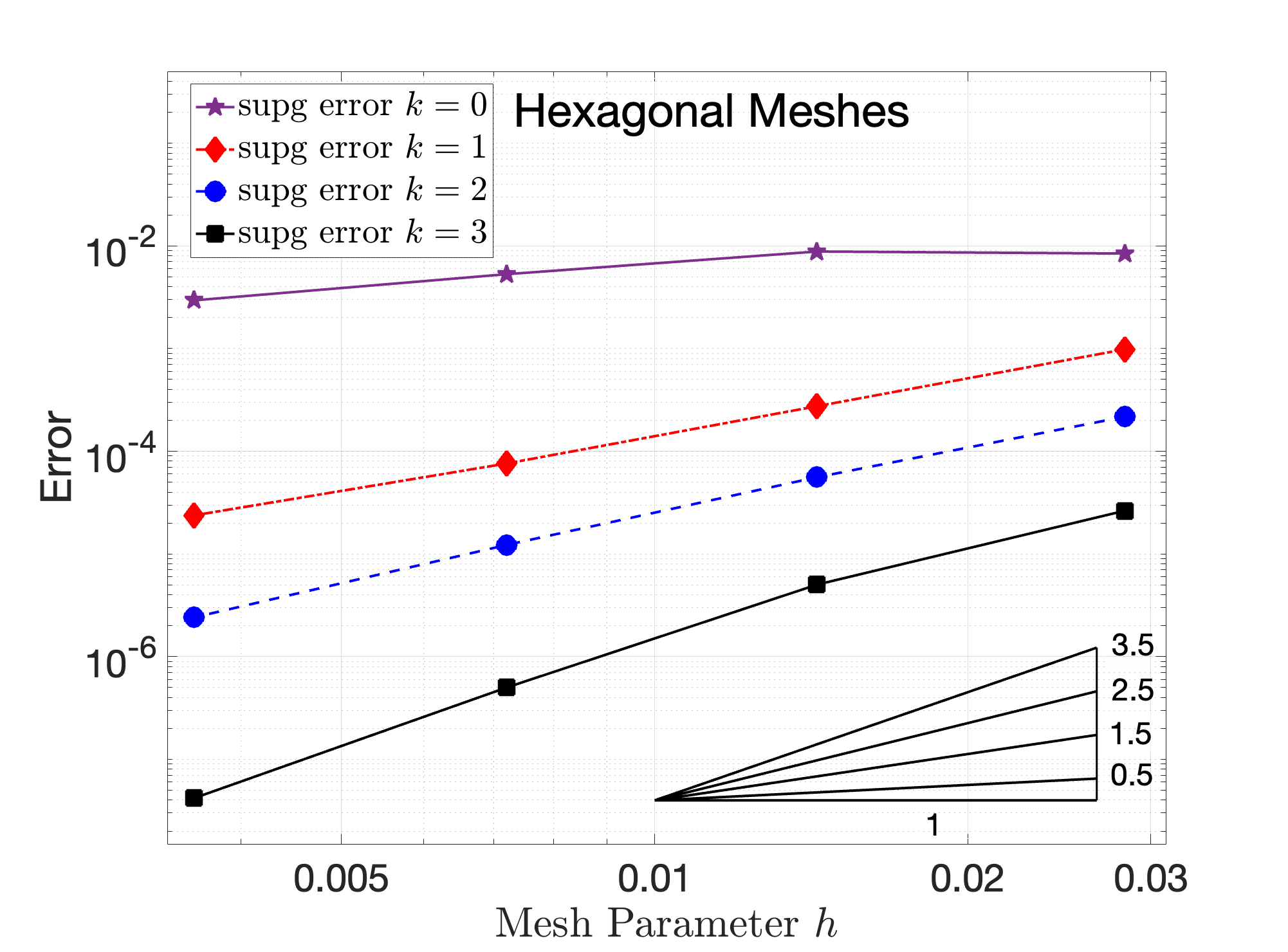}}
		\caption{Convergence histories for the error (left part - (a), (c) and (e)) in $\tnorm{\cdot}_{\rm LP}$ norm of Example~\ref{test:smooth} on the triangular, Cartesian, and hexagonal meshes for $\epsilon=10^{-8}$ and in the supg norm (right part - (b), (d) and (f)) $\|\cdot\|_{\rm supg}$ for $k=0,1,2,3$.}
		\label{fig:Smooth_Conv_His_Triag_Cart_Hexa}
	\end{center}
\end{figure}

\end{example}

\begin{example}[Boundary layer problem]\label{test:exp_layer}
Consider the test problem \eqref{Oseen_eqn_intro} for $\epsilon=10^{-2}$ with $\b=(1,1)$ and $\sigma=1$. Let the exact solution of velocity and pressure be given by  
\begin{align*}
u(x,y)=\begin{bmatrix}
2x^2y(e^{\lambda(x-1)} - 1)^2(e^{\lambda(y-1)} - 1)^2 + 2x^2y^2\lambda e^{z(y-1)}(e^{z(x-1)} - 1)^2(e^{\lambda(y-1)} - 1)\\
2xy^2(e^{\lambda(x-1)} - 1)^2(e^{\lambda(y-1)} - 1)^2 + 2x^2y^2\lambda e^{\lambda(x-1)}(e^{\lambda(x-1)} - 1)(e^{\lambda(y-1)} - 1)^2
\end{bmatrix},
\end{align*}
\begin{align*}
p(x,y)=e^{x+y}-(e-1)^2,
\end{align*}
where $\lambda=\frac{1}{2\sqrt{\epsilon}}$.

With the above problem setups, we perform numerical tests on the triangular, Cartesian, and hexagonal mesh families as described in the previous Example~\ref{test:smooth}. In the left part of Figure~\ref{fig:Exp_Layer_Conv_His_Triag_Cart_Hexa} ((a), (c) and (e)), we plot the error $\tnorm{(\underline{\boldsymbol{e}}_h^{\u}, e_h^{p})}_{\rm LP}$ as a function of meshsize $h$ for polynomial degrees $k=1,2,3$ with $\epsilon=10^{-2}$ on triangular, Cartesian, and hexagonal meshes. In the right part of Figure~\ref{fig:Exp_Layer_Conv_His_Triag_Cart_Hexa} ((b), (d) and (f)), we plot the supg error $\tnorm{(\underline{\boldsymbol{e}}_h^{\u}, e_h^{p})}_{\rm supg}$ as a function of meshsize $h$ for polynomial degrees $k=1,2,3$ with $\epsilon=10^{-2}$ on triangular, Cartesian and hexagonal meshes. In both cases, the convergence rates are approximately $h^{k+1/2}$ for $\epsilon=10^{-2}$.
\begin{figure}
	\begin{center}
		\subfloat[]{\includegraphics[height=0.4\textwidth,width=0.5\textwidth]{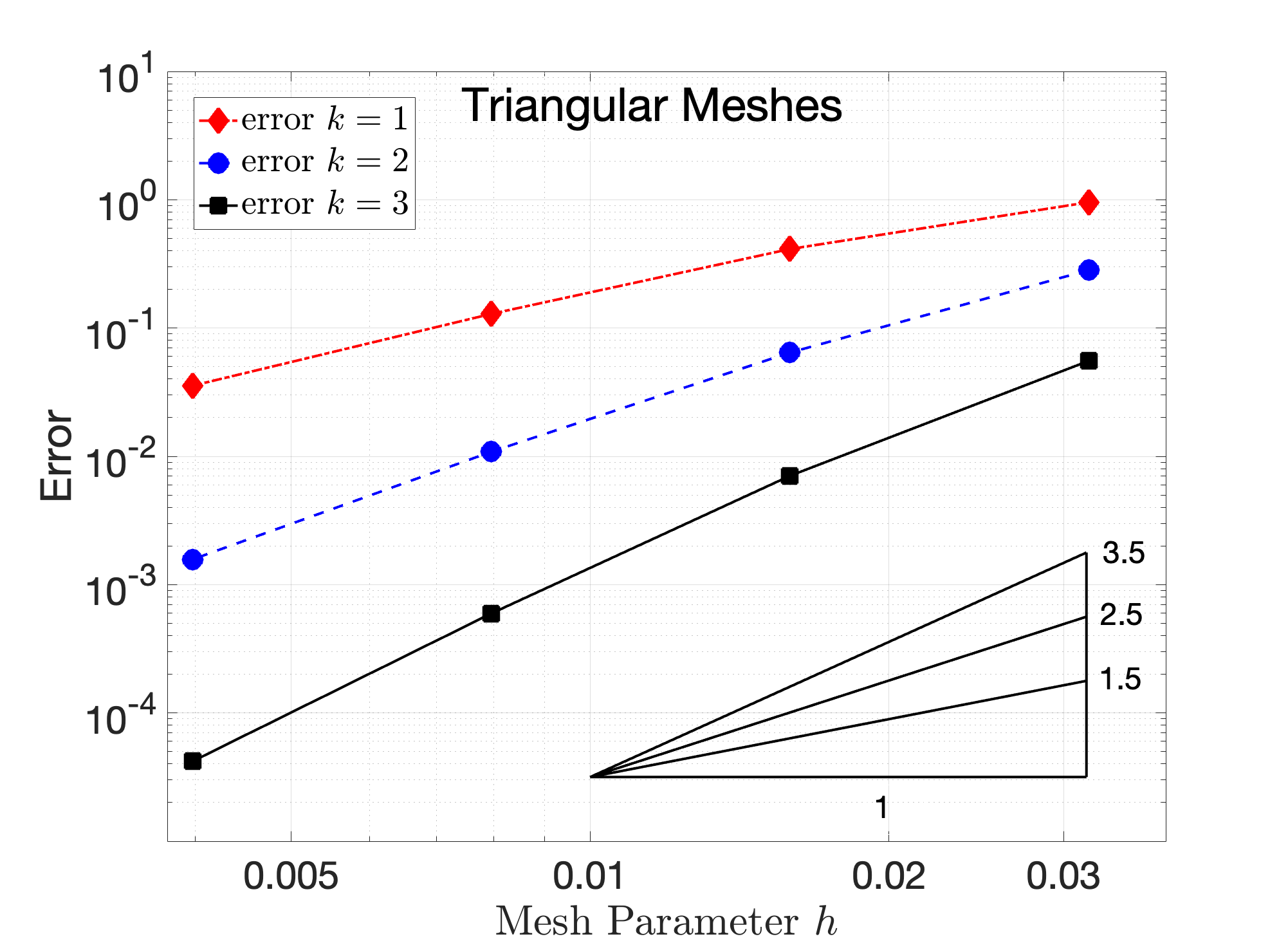}}
  		\subfloat[]{\includegraphics[height=0.4\textwidth,width=0.5\textwidth]{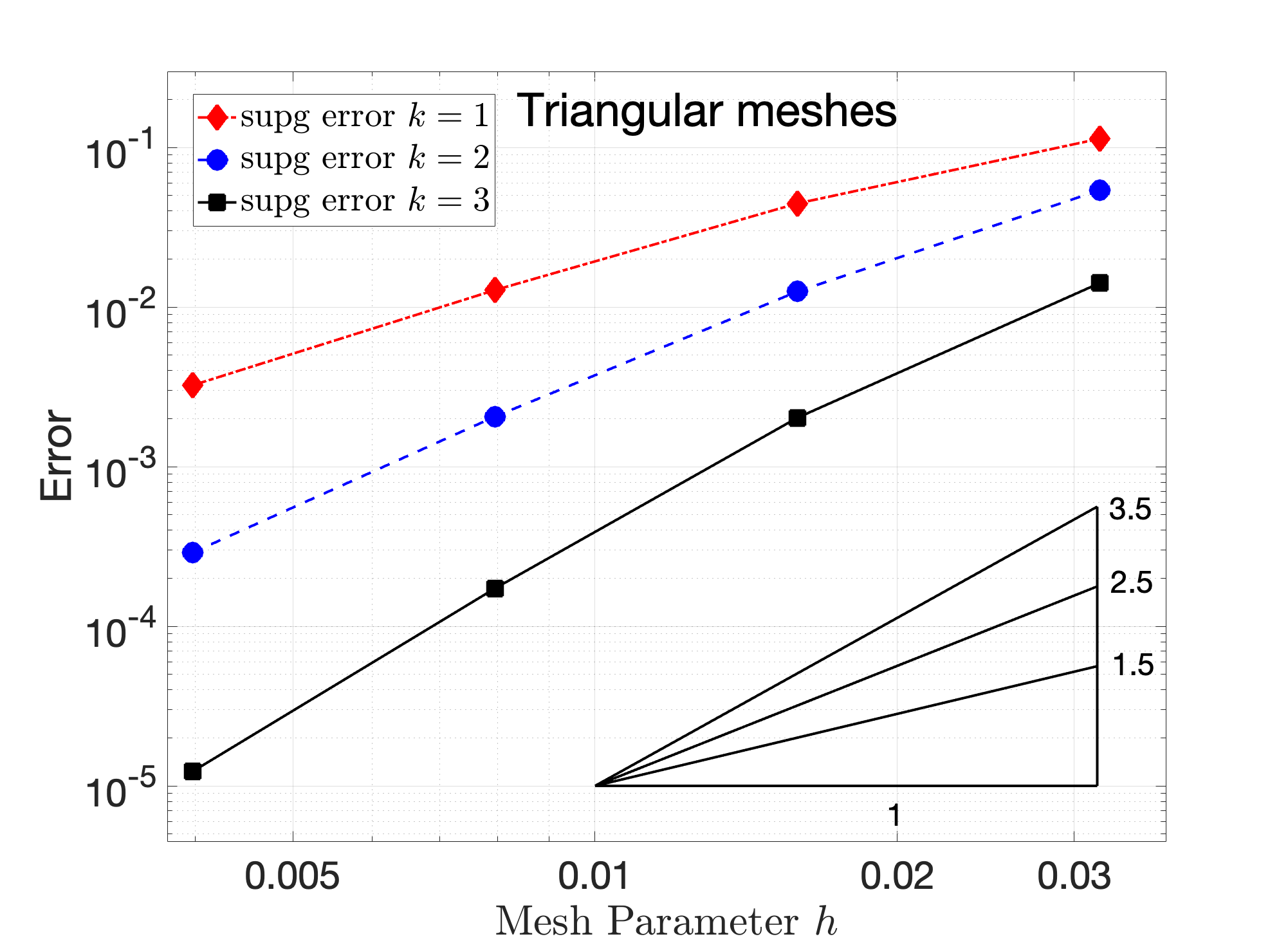}}\\
		\subfloat[]{\includegraphics[height=0.4\textwidth,width=0.5\textwidth]{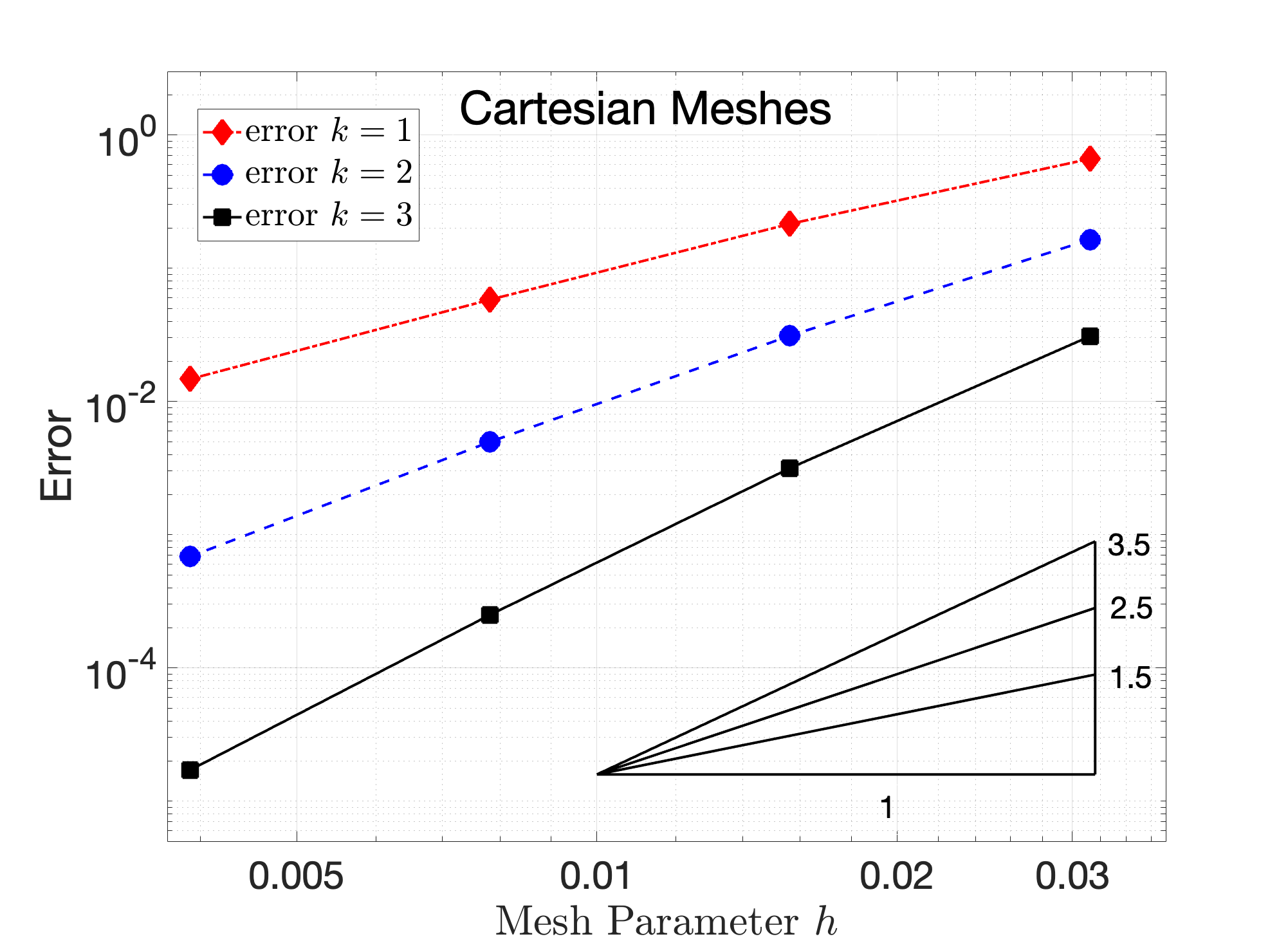}}
  		\subfloat[]{\includegraphics[height=0.4\textwidth,width=0.5\textwidth]{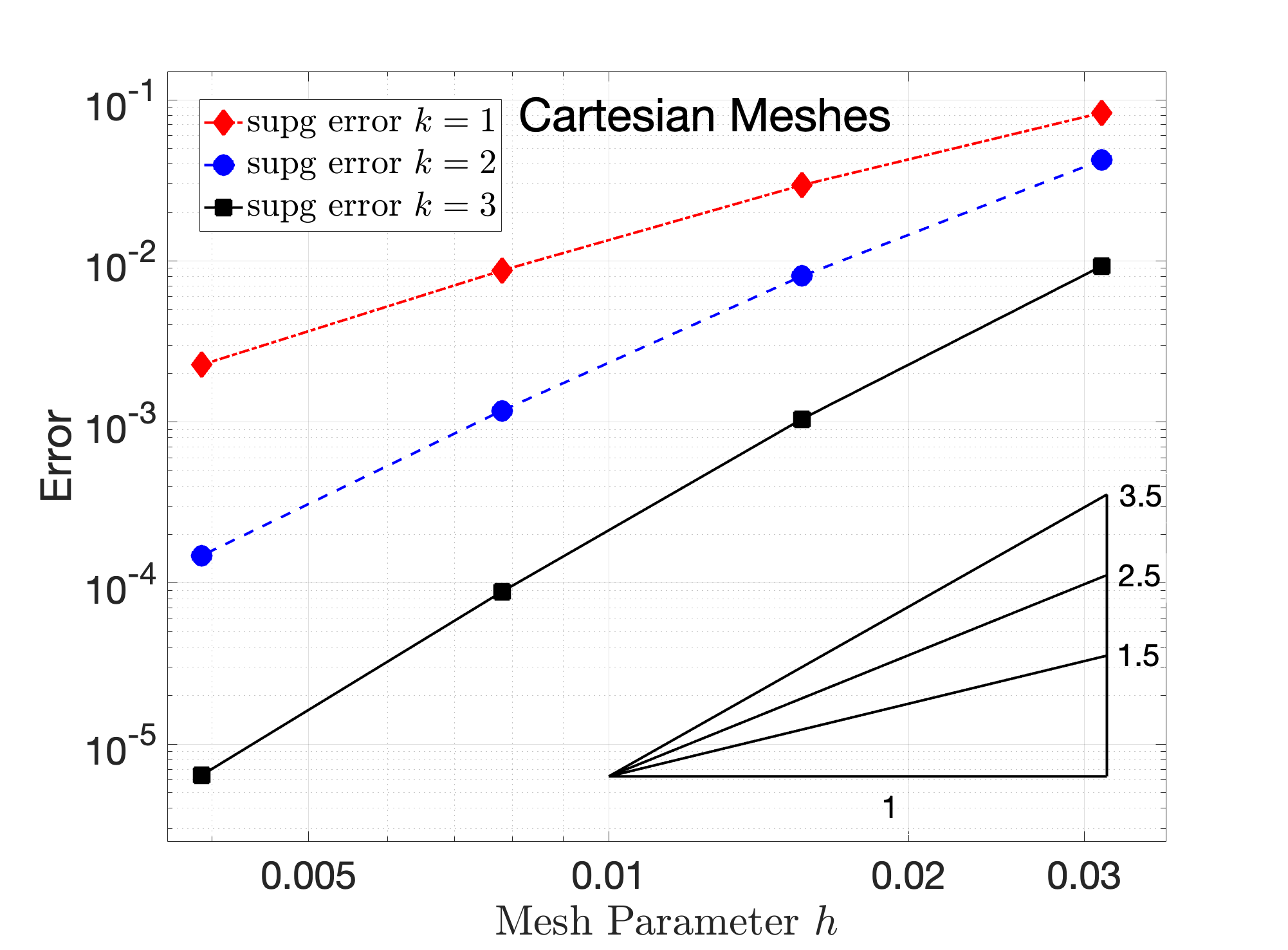}}\\
		\subfloat[]{\includegraphics[height=0.4\textwidth,width=0.5\textwidth]{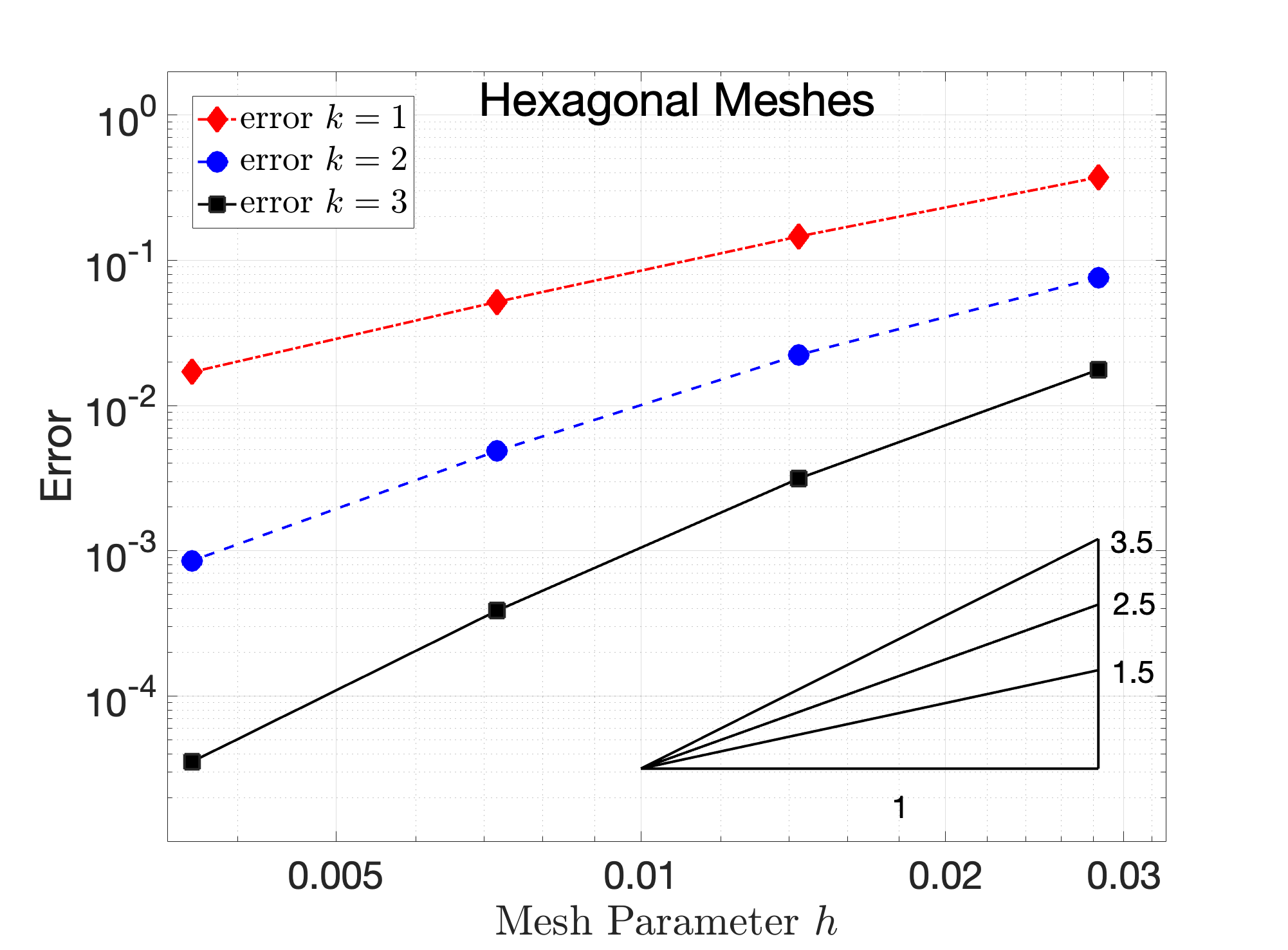}}
  		\subfloat[]{\includegraphics[height=0.4\textwidth,width=0.5\textwidth]{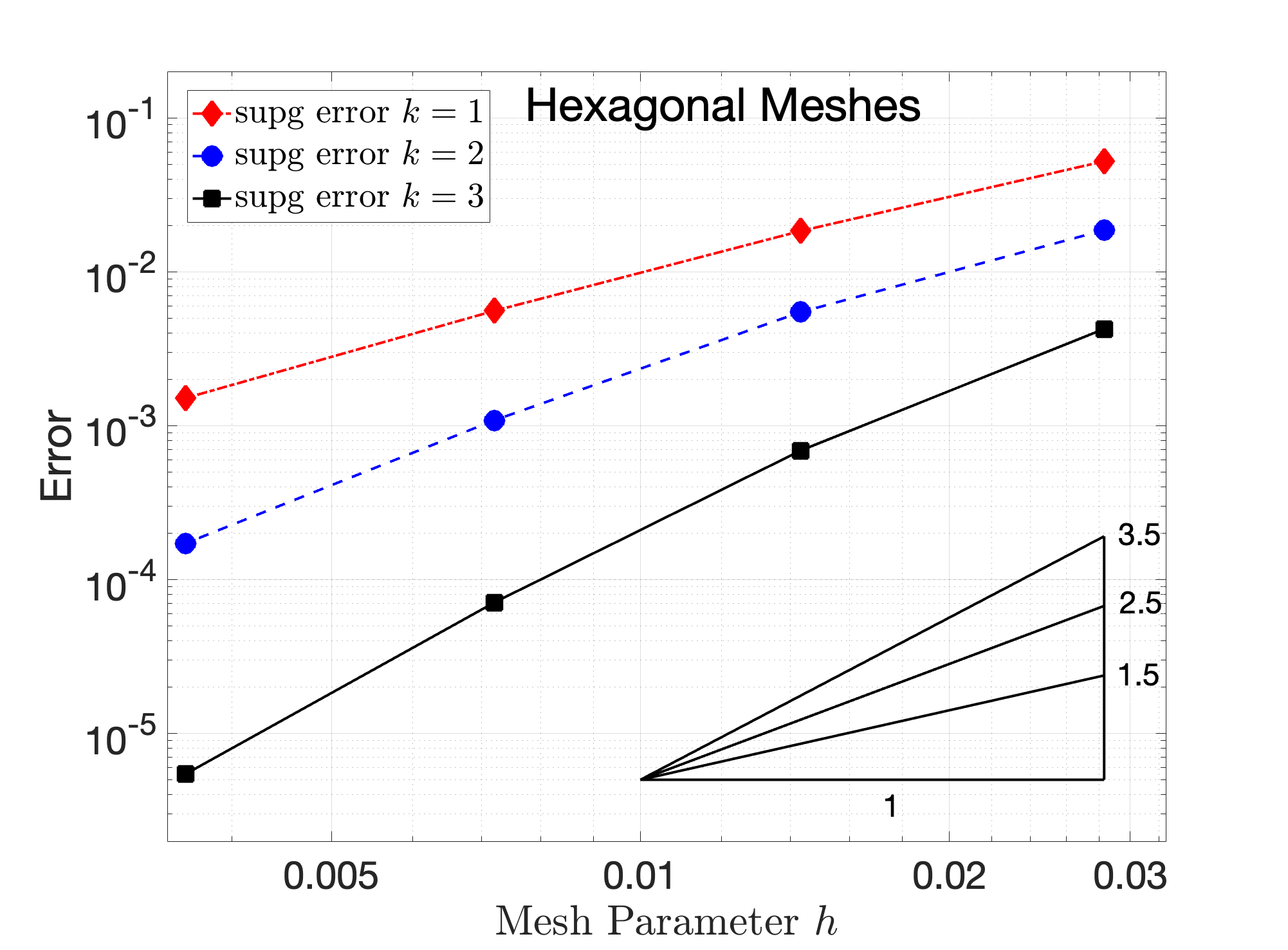}}
		\caption{Convergence histories for the error (left part - (a), (c) and (e)) in $\tnorm{\cdot}_{\rm LP}$ norm of Example~~\ref{test:exp_layer} on the triangular, Cartesian, and hexagonal meshes for $\epsilon=10^{-2}$ and in the supg norm (right part - (b), (d) and (f)) $\|\cdot\|_{\rm supg}$ for $k=1,2,3$.}
		\label{fig:Exp_Layer_Conv_His_Triag_Cart_Hexa}
	\end{center}
\end{figure}
\end{example}

\section{Conclusions}
In this article, we have considered a local projection stabilised Hybrid High-Order method for the Oseen problem. We have proved a strong stability result for our discrete formulation under the SUPG norm. The error analysis shows the optimal order of convergence $(k+1/2)$, assuming that velocity and pressure have regularity $H^{k+1}$ globally and $H^{k+2}$ locally. In the second example with very small $\epsilon$, we obtain suboptimal convergence rates for both the HHO-LPS method and the HHO method discussed in \cite{Oseen_HHO}. It is our presumption that the issue of the suboptimal convergence rate in example \ref{test:exp_layer} with very small $\epsilon$ can be resolved using an adaptive approach. This will be considered in future work.

\bigskip

\noindent\textbf{Acknowledgements}\\
The first author acknowledges the DST-SERB-MATRICS grant MTR/2023/000681 for financial support and DST FIST grant SR/FST/MS II/2023/139 for partial financial support.

\medskip

\section{Appendix}
Note that the global divergence reconstruction operator $D_h^k$ defined as $D_h^k\rvert_T=D_T^k$ produces globally  $L^2_0(\Omega)$ function. This can be checked simply by taking $p=1$ on each cell $T \in \cT_h$ and using the fact that $\bv_F=0$ on the boundary edges for a hybrid function $\bvl_h=(\bvl_T,\bvl_F) \in \Ul_{h,0}^k$.
\begin{align*}
\sum_{T \in \cT_h}(D_T^k\bvl_T,1)_T=\sum_{T \in \cT_h}\sum_{F \in \cF_T}(\bv_F \cdot \bn_{TF},1)_F=0.
\end{align*}
Along with this, $D_h^k$ is also a continuous linear surjective operator from $\Ul_{h,0}^k$ to $P_h^k$. Moreover, we have the following lemma.

\begin{lemma}\label{velocity_pressure_infsup}
There exists $\alpha > 0$ independent of the meshsize such that for all $p_h \in P_h^k$ we have the following
\begin{align}\label{velocity_pressure_infsup_eqn}
\sup_{\bvl_h \in \Ul_{0,h}^k}\frac{B_h(\bvl_h,p_h)}{\norm{\bvl_h}_{1,h}} \geq \alpha \norm{p_h}.
\end{align}
\end{lemma}
\begin{proof}
It is well known that the velocity/pressure pair $[H^1_0(\Omega)]^d/L^2_0(\Omega)$ is stable and the continuous divergence operator is surjective. Therefore, for any $p_h \in P_h^k$ one can have a unique $\bv \in \V$ such that $\norm{\bv}_1\leq C\norm{p_h}$ and $\div \bv =-p_h$. Using this and the relation $\pi_T^{k}(\div \bv)=D_T^k(\Il_h(\bv))$ we get
\begin{align}\label{velocity_pressure_infsup_eqn_1}
\norm{p_h}^2_{\Omega} = (-\div \bv,p_h)_{\Omega}= B_h(\Il_h^k(\bv),p_h)
\end{align}
Moreover, using the approximation properties of $\pi_T^k$ and $\pi_F^k$ we have
\begin{align}\label{velocity_pressure_infsup_eqn_2}
\norm{\Il_h(\bv)}^2_{1,h} &=\sum_{T \in \cT_h}\Big(\norm{\nabla \Il_T(\bv)}^2_T+\sum_{F \in \cF_T} \int_{F}\frac{1}{h_F} (\pi_F^{k}(\bv)-\pi_T^{k}(\bv))^2\Big) \nonumber\\
&\leq C \norm{\bv}_1^2\leq C\norm{p_h}^2.
\end{align}
Using the relations \eqref{velocity_pressure_infsup_eqn_1} and \eqref{velocity_pressure_infsup_eqn_2} we arrive at \eqref{velocity_pressure_infsup_eqn} with $\alpha = 1/C$.
\end{proof}
Since $P_h^k$ is a reflexive space and the discrete divergence operator $D_h$ is surjective, the converse of the Lemma A.42 in \cite{Ern:2004:FEMbook} says that for any $p_h \in P_h^k$ there exists $\bvl_h \in \Ul_{0,h}^k$ such that 
\begin{align}\label{eqdivp}
D_h^k(\bvl_h)=p_h  \quad \text{and} \quad \norm{\bvl_h}_{1,h}^2 \leq C\norm{p_h}^2, ~ \text{C depends on } \Omega.
\end{align}

\bibliographystyle{siam}
\bibliography{HHO,Oseen}

\end{document}